\documentclass[11pt, oneside]{article} 
\makeatletter
\def\cl@chapter{\@elt {theorem}}
\makeatother
\usepackage{graphicx}
\usepackage{mathptmx}      % use Times fonts if available on your TeX system
\usepackage{geometry}                		% See geometry.pdf to learn the layout options. There are lots.
\usepackage[mathlines]{lineno}
\usepackage{listings}                 		% ... or a4paper or a5paper or ... 
\usepackage{graphicx}				% Use pdf, png, jpg, or eps§ with pdflatex;		
\usepackage[misc]{ifsym}
\usepackage{color}
\usepackage{amssymb}
\usepackage{amsmath}
\usepackage{amsthm}
\usepackage{mathtools}
\usepackage{xfrac}
\usepackage{bm}
\usepackage{aliases}
\usepackage[utf8]{inputenc}
\usepackage[english]{babel}
\usepackage{textcomp}
\usepackage{stmaryrd}
\usepackage{todonotes}
\usepackage{siunitx}[=v2]
\usepackage{commath}
\usepackage[affil-it]{authblk}
\usepackage[outdir=./Figures/epstopdf/]{epstopdf}
\numberwithin{equation}{section}
\usepackage{siunitx}
\usepackage{csquotes}
\usepackage[unicode, pdfpagelabels,bookmarks,hyperindex,hyperfigures]{hyperref}
\usepackage{cleveref}
\usepackage[outdir=./Figures/epstopdf/]{epstopdf}
\numberwithin{equation}{section}
\usepackage{siunitx}
\graphicspath{{Figures/},{Figures/eps/},{{Figures/jpeg/}}}
\usepackage[numbers]{natbib}
\author[1]{Arijit Hazra}
\author[2]{Dinshaw S. Balsara}
\author[3]{Praveen Chandrashekar}
\author[4]{Sudip K. Garain}

\affil[1]{\footnotesize Department of Mechanical Engineering, Indian Institute of Technology, Palakkad, India Email: \texttt{ahazra@iitpkd.ac.in}, \texttt{mailtohazra@gmail.com}}
\affil[2]{\footnotesize Department of Physics, University of Notre Dame, USA. Email: \texttt{dbalsara@nd.edu}}
\affil[3]{\footnotesize TIFR Center for Applicable Mathematics, Bangalore, India. Email:\texttt{praveen@math.tifrbng.res.in}}
\affil[4]{\footnotesize   Department of Physical Sciences
	Indian Institute of Science Education and Research Kolkata, India. Email: \texttt{sgarain@iiserkol.ac.in}}

\title{Multidimensional Generalized Riemann Problem Solver for Maxwell’s Equations}
\usepackage[unicode, pdfpagelabels,bookmarks,hyperindex,hyperfigures]{hyperref}
\usepackage{cleveref}
\begin{document}
%\linenumbers
\maketitle
\begin{abstract}
	
	Approximate multidimensional Riemann solvers are essential building blocks in designing globally constraint-preserving finite volume time domain (FVTD) and discontinuous Galerkin time domain (DGTD) schemes for computational electrodynamics (CED).  In those schemes, we can achieve high-order temporal accuracy with the help of Runge-Kutta or ADER time-stepping. This paper presents the design of a multidimensional approximate Generalized Riemann Problem (GRP) solver for the first time. The multidimensional Riemann solver accepts as its inputs the four states surrounding an edge on a structured mesh, and its output consists of a resolved state and its associated fluxes. In contrast, the multidimensional GRP solver accepts as its inputs the four states and their gradients in all directions; its output consists of the resolved state and its corresponding fluxes and the gradients of the resolved state. The gradients can then be used to extend the solution in time. As a result, we achieve second-order temporal accuracy in a single step.
	
	In this work, the formulation is optimized for linear hyperbolic systems with stiff, linear source terms because such a formulation will find maximal use in CED. Our formulation produces an overall constraint-preserving time-stepping strategy based on the GRP that is provably L-stable in the presence of stiff source terms. We present several stringent test problems, showing that the multidimensional GRP solver for CED meets its design accuracy and performs stably with optimal time steps. The test problems include cases with high conductivity, showing that the beneficial L-stability is indeed realized in practical applications.
	
	{\bf Keywords}: {Conservation laws, Hyperbolic partial differential equations, Multidimensional Riemann problem, Maxwell's equations}
\end{abstract}
	
%%%---------------------------------Introduction --------------------------------------------------------------%%%%
\section{Introduction}
\label{sec:intro}
Computational Electrodynamics (CED) which deals with the numerical solution of Maxwell’s equations, plays a vital role in many problems in science and engineering. The finite-difference time-domain (FDTD) method \cite{yee1966numerical,taflove1975numerical,taflove1988review,taflove2000computational,taflove1999finite} has been a mainstay of CED applications. The primary strength of FDTD stems from its use of a beneficial staggering of the electric and magnetic fields to ensure that the global constraints (inherent in Gauss’s law and the absence of magnetic monopoles) are discretely represented on the computational mesh. FDTD is globally constraint-preserving. However, the {{primary}} weakness of {{standard }} FDTD stems from the fact that it is restricted to second-order accuracy, especially when electromagnetic radiation interacts with material media. 

The differential form of Maxwell’s equations has a dissipationless and dispersionless limit. As a result, it is beneficial for numerical schemes to be as dissipationless and dispersionless as possible. This has given rise to the discontinuous Galerkin time domain (DGTD) methods  \cite{angulo2015discontinuous,hesthaven2002nodal,hesthaven2007nodal,chen2012discontinuous,chen2009non,ren2015new,wang2017new,sun2017novel}.
Such methods do not satisfy the constraints in a global sense; though some of them do satisfy the constraints locally within each element. Even so, since they are based on discontinuous Galerkin methods, their strong point is that they can reach high orders of accuracy. It is very desirable to retain good traits of the FDTD and DGTD schemes discussed above.

In an effort to design CED schemes that offer the best of both worlds – global constraint preservation from FDTD and higher order from DGTD – we have embarked on an effort to design such schemes. Therefore, finite volume time domain (FVTD) schemes that globally preserve constraints and also attain high order of accuracy were presented in \cite{balsara_high-order_2016,balsara_computational_2017,balsara_computational_2018}.
DGTD schemes with those same favorable attributes were presented in  \cite{balsara2017neumann,balsara2019neumann,hazra_globally_2019}.
The two central ingredients of those schemes are a high order constraint-preserving reconstruction of vector fields  \cite{balsara_divergence-free_2001,balsara_secondorderaccurate_2004,balsara_divergence-free_2009,balsara_divergence-free_2015,balsara_high-order_2016,balsara_computational_2017,balsara_computational_2018,xu_divergence-free_2016}
and multidimensional Riemann solvers  \cite{balsara_multidimensional_2010,balsara_two-dimensional_2012,balsara_multidimensional_2014-1,balsara2015three,balsara2014multidimensional,balsara_multidimensional_2015,balsara2016two,balsara2017multidimensional}
The constraint-preserving reconstruction provides spatially high order accuracy. The multidimensional Riemann solver folds in the essential physics that electromagnetic phenomena are mediated by wave propagation that invariably occurs in all directions. It also gives us a natural, physics-based approach for obtaining the electric and magnetic fields at the edges of the computational mesh. 
%The high order spatial accuracy should be matched by high order temporal accuracy.

{{Furthermore, Maxwell's equations have symplectic and multi-symplectic structures. Considering this, Leapfrog time integration has been the chosen strategy for standard FDTD as it is a form of symplectic integrator. However, it is well-known that FDTD results in high level of dispersion error \cite{taflove2000computational}. As this numerical dispersion error accumulates over time, simulation of long-term behaviour and long-duration electromagnetic wave propagation with FDTD requires an extremely fine mesh, and finer mesh in conjunction with courant stability criteria results in prohibitively high computational time for such simulations with FDTD. There have been several efforts to reduce the dispersion error by modifying FDTD \cite{smith_hierarchy_2012}. However, considering all the desirable features of higher order numerical methods for CED, such as higher order spatial and temporal accuracy, ability to handle complex geometry, low dispersion error, higher-order CED schemes generally tend to use low-storage five-stage fourth-order Explicit Runge-Kutta method (LSERK4)\cite{williamson_low-storage_1980,hu_low-dissipation_1996,berland_high-order_2007,hesthaven2007nodal,niegemann_efficient_2012,diehl2010comparison}, strongly stability-preserving Runge-Kutta (SSPRK) \cite{sarmany_dispersion_2007,chen_high-order_2005,hazra_globally_2019,balsara2019neumann} or ADER (Arbitrary DERivatives in space and time) \cite{dumbser_unified_2008,dumbser_ader-weno_2013,balsara_efficient_2009,balsara_efficient_2013} time-discretizations \cite{taube_high-order_2009,balsara_computational_2018}.
		
		To compare briefly the computational complexities of Runge-Kutta and ADER time integration schemes for CED,}~we first note that each stage in a Runge-Kutta time-discretization is only first order accurate in time.  For higher-order constraint-preserving time evolution of CED with Runge-Kutta schemes is, therefore, obtained by the application of a multidimensional Riemann solver at the edges of the mesh in order to obtain the edge-collocated integrals of the electric and magnetic fields. Thus each stage of a Runge-Kutta time-discretization is relatively inexpensive, but the overall scheme can be more expensive because multiple stages are used. Since the CED equations can have stiff source terms, the inclusion of stiff source terms can also add to the cost of a Runge-Kutta time-discretization. The ADER update only requires a single stage ADER formulation within each zone to make an “in-the-small” evolution within each zone. Once this is available, constraint-preserving time evolution of CED can be obtained with volumetrically-based ADER schemes by invoking a multidimensional Riemann solver at the edges of the mesh. However, for volumetrically-based ADER schemes the space-time ADER construction within a zone can itself be quite expensive. The treatment of stiff source terms also adds to the cost of an ADER scheme.
	
	The utility of a GRP approach stems from the fact that a GRP can utilize not just the input states, but also their gradients. Realize that those gradients are always available, and they can always be provided by the spatial reconstruction. The intricacy in designing a GRP solver consists of finding ways to take the gradients of the input states and using them to obtain gradients in the resolved state. Once the gradients in the resolved state are obtained, one can obtain at least a second order accurate update in one stage. While a few exact and approximate GRP solvers have been designed that go beyond second order accuracy \cite{le_floch_asymptotic_1988,titarev2002ader,titarev2005ader,toro2001towards,toro_solution_2002,qian_generalized_2014,montecinos2014reformulations,wu2014third,goetz_novel_2016,goetz_family_2018}, 
	the majority of GRP solvers have been restricted to second order in time \cite{ben-artzi_second-order_1984,ben-artzi_upwind_1986,ben-artzi_generalized_2003,ben-artzi_generalized_1989,ben1990computation,bourgeade_asymptotic_1989,balsara_efficient_2018}. All the GRP constructions that we know of have been one-dimensional. Since multidimensional Riemann solvers have begun to play such an important role in CED, and also other fields, it is of great interest to obtain generalized Riemann problem versions of the same.
	
	In globally-constraint preserving schemes for CED, we apply a multidimensional Riemann solver to the edges of the mesh. Such multidimensional Riemann solvers have been designed \cite{balsara_multidimensional_2010,balsara_two-dimensional_2012,balsara_multidimensional_2014-1,balsara2015three,balsara_multidimensional_2014,balsara2016two,balsara2017multidimensional}.
	However, as far as we know, this is the first effort to formulate a multidimensional generalized Riemann problem solver that works seamlessly. The goal of this first paper is to design a multidimensional generalized Riemann solver for CED. We choose CED because it is a linear hyperbolic system and it is very beneficial to study the problem in the context of a linear system before tackling the fully non-linear case. The fully non-linear case will be formulated in a subsequent paper. We formulate the problem so that it can be used for any general linear hyperbolic system, but we also specialize our results for CED.
	
	CED, just like aeroacoustics, is very special in that most applications are linear. If non-linearities are present, they are usually mild. But that only changes the emphasis of the solution methodology. Because waves can propagate without dissipation or dispersion in electrodynamics and aeroacoustics, a substantial premium is placed on minimizing numerical dispersion and dissipation.
	There has been a growing realization that the availability of GRP solvers can lead to a new generation of {{low-dissipation, low-dispersion}}~Taylor Series-based (TS-based) schemes \cite{li_two-stage_2016,christlieb_explicit_2016,grant_strong_2019}, though that field is perhaps still emergent. The schemes are referred to as Taylor series-based because the GRP solver delivers not just the numerical flux but also its derivative in time. 
	%Multi-stage TS-based schemes that minimize dissipation and dispersion while retaining high order of temporal accuracy have not yet been developed. For that reason, we restrict this paper to second order of accuracy in space and time. 
	{{The novelty of our work lies in presenting a multidimensional GRP solver, which can be an essential building block for the development of low dissipation, low dispersion TS-based schemes for CED, aeroacoustics and other analogous fields. We also show how linear stiff source terms can be included in the multidimensional GRP solver.}}
	
	The rest of the paper is organized as follows. In \Cref{sec:maxwellEqn} we describe Maxwell’s equations and globally constraint-preserving solution methods for those equations. The multidimensional GRP solver is described in \Cref{sec:multidGRP}. \Cref{sec:Implementation} gives a pointwise strategy for implementation. \Cref{sec:accuracy} provides accuracy analysis; \Cref{sec:testcases} provides several stringent test problems. \Cref{sec:conclusion} draws some conclusions.
		
	%%%--------------------------------- Maxwell's equations --------------------------------------------------------------%%%%
	
	\section{Maxwell equation} \label{sec:maxwellEqn}
	We split this Section into two parts. \Cref{subsec:intromaxwell} introduces Maxwell’s equations. \Cref{subsec:globalsolMaxwell} describes their globally constraint-preserving numerical solution using a GRP solver.
	
	\subsection{Introduction to Maxwell’s equations} \label{subsec:intromaxwell}
	
	The equations of CED can be written as two evolutionary curl-type equations for the magnetic induction and the electric displacement. The first of these is Faraday’s law, given by,
	%\begin{frame}{Maxwell's equations}
	\begin{linenomath*}
		\begin{align}
			\dpd{\bm{B}}{t} + \nabla \times \bm{E} = -\M,
			\label{eqn:faradayslaw}
		\end{align}
	\end{linenomath*}
	where $\B$ is the magnetic induction (or magnetic flux density),  $\E$ is the electric field and $\M$ is the the magnetic current density. The magnetic current density is zero for any physical domain. The second evolutionary equation for the electric displacement is the extended Ampere’s law, given by
	\begin{align}
		\dpd{\bm{D}}{t} - \nabla       \times \bm{H} = -\bm{J},
		\label{eqn:amperslaw}
	\end{align}
	where $\D$  is the electric displacement (or electric flux density),  $\bH$ is the magnetic field vector and  $\J$ is the electric current density. The structure of the above two equations is such that the magnetic induction and the electric displacement also satisfy the following two non-evolutionary involution constraints, given by
	\begin{align}
		\nabla \cdot \bm{B}= \rho_M,
	\end{align}
	and
	\begin{align}
		\nabla \cdot \D= \rho_E.
	\end{align}
	Here $\rho_M $ and $\rho_E$ are the magnetic and electric charge densities. For any physical medium  $\rho_M=0$ since magnetic monopoles do not exist.
	
	The involutionary nature of the above equations ensures that the electric charge density satisfies the equation
	\begin{align}
		\dpd{\rho_E}{t} + \nabla \cdot \bm{J} = 0,
	\end{align}
	and the magnetic charge density satisfies the equation
	\begin{align}
		\dpd{\rho_M}{t} + \nabla \cdot \bm{M} = 0.
	\end{align}
	In material media we also have the constitutive relations 
	\begin{align}
		\bm{B} = \mu \bm{H},
	\end{align}
	and
	\begin{align}
		\bm{D} = \varepsilon \bm{E},
	\end{align}
	where $\mu$ is a \num{3 x 3} permeability tensor and $\varepsilon$ is the analogous \num{3 x 3} permittivity tensor. For most material media, these tensors are diagonal. The eigenstructure of the hyperbolic system is most easily found for the diagonal case, where we make the simplifying assumption  $\varepsilon = \mathrm{diag}\{\varepsilon_{xx}, \varepsilon_{yy}, \varepsilon_{zz}\}$ and $\mu = \mathrm{diag}\{\mu_{xx}, \mu_{yy}, \mu_{zz}\}$. The corresponding eigenstructure has been catalogued in Sub-section II.2  of \cite{balsara_computational_2017}. We will also need the inverses of the permittivity and permeability tensors. These \num{3 x 3}  inverse matrices will also be symmetric, and we denote them as  $\ieps$ and $\imu$ . 
	
	The current density is related to the electric field via
	\begin{align}
		\bm{J} = \sigma \bm{E},
	\end{align}
	where $\sigma$ is the conductivity.  Similarly, the magnetic current density is related to the magnetic field via
	\begin{align}
		\M = \sigma^* \bH,
	\end{align}
	where  is the equivalent magnetic loss, which is again zero in physical media, but may be non-zero when imposing boundary conditions in CED.
	
	\subsection{Globally constraint-preserving numerical solution of Maxwell’s equations}
	\label{subsec:globalsolMaxwell}
	\begin{figure}
		\begin{center}
			\includegraphics[width=0.80\textwidth]{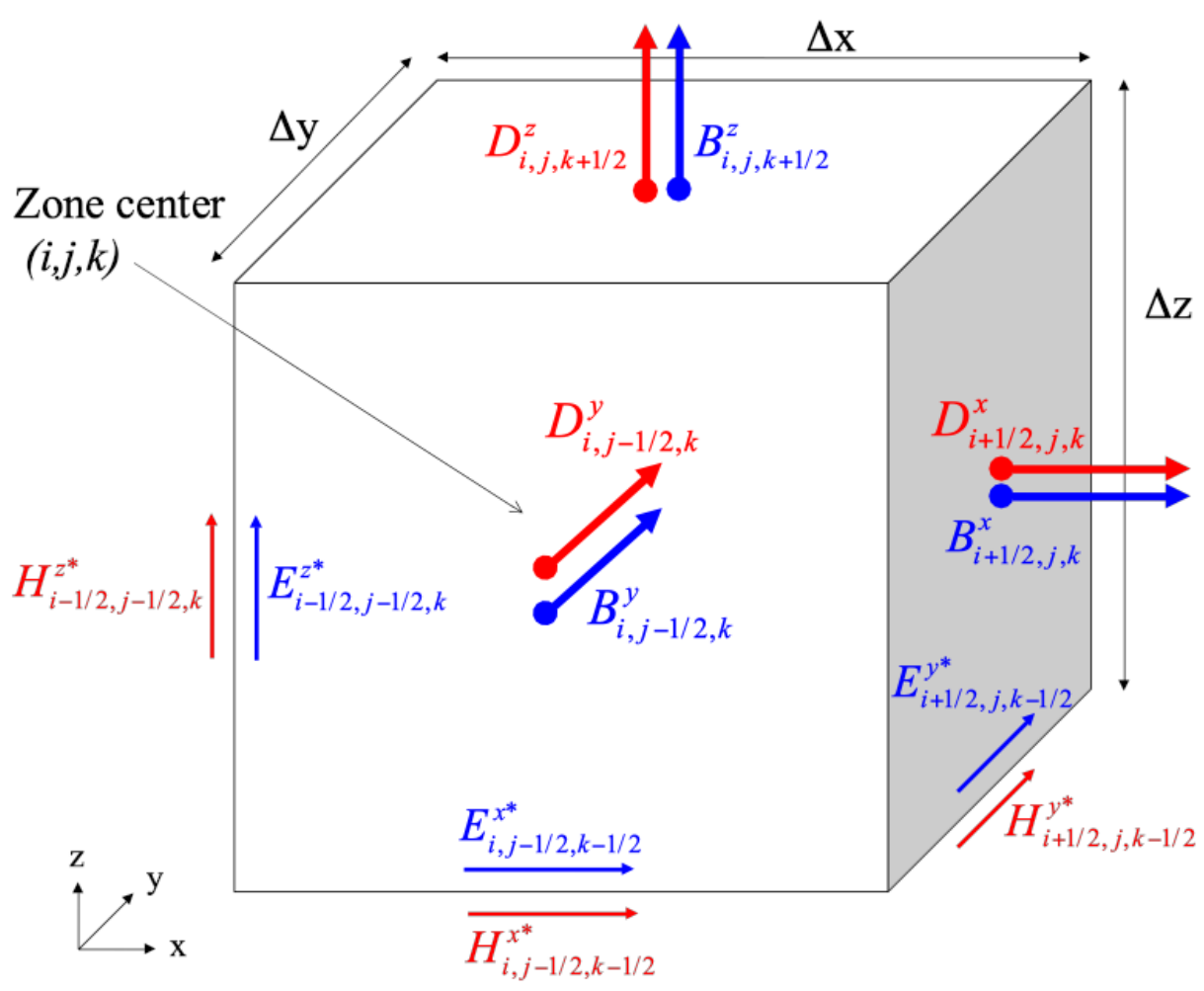} 
			\caption{Schematic diagram depicting the collocation of the primal and dual variables of Maxwell's equations. Primal variables of the scheme, given by the normal components of the magnetic induction and electric field displacement, are facially-collocated. They undergo an update from Faraday’s law and the generalized Ampere’s law, respectively. The components of the primal magnetic induction vector are shown by the thick blue arrows, while the components of the primal electric displacement vector are shown by the thick red arrows. The edge-collocated electric fields, which are used for updating the facial magnetic induction components, are shown by the thin blue arrows close to the appropriate edge. The edge-collocated magnetic fields, which are used for updating the facial electric displacement components, are shown by the thin red arrows close to the appropriate edge. }
			\label{fig:3ddivfree}
		\end{center}
	\end{figure}
	The facially-collocated normal components of the electric displacement and the magnetic induction constitute the primal variables of our scheme 
	In \Cref{fig:3ddivfree}, these vector fields are shown by the thick red arrow and the thick blue arrow, respectively, in each of the faces of the cuboidal element. In a finite-volume sense, these primal variables are actually taken to be facial averages of the normal components of the electric displacement and the magnetic induction. The overall task consists of finding the edge-collocated components of the magnetic field vector and the electric field vector, shown in \Cref{fig:3ddivfree}. These are shown with the thinner red arrow and the thinner blue arrow, respectively, next to the edges of the zone shown in \Cref{fig:3ddivfree}. In a finite-volume sense, these are actually averages in one space dimension (taken to be the length of the element’s edge) and the time dimension (evaluated over the timestep). The order of spatial reconstruction of the electric displacement and the magnetic induction then determines the order of spatial accuracy of our numerical scheme. At second order, volumetric ADER schemes of the sort designed in \cite{dumbser_unified_2008} and \cite{balsara_efficient_2009} can indeed provide a one-step update. However, a similar one-step update can be obtained using the multidimensional GRP solver designed here. 
	
	A single step constraint-preserving update for the entire set of CED equations, consistent with the curl-type update in \Cref{eqn:faradayslaw,eqn:amperslaw}  can be written at each face of the zone shown in \Cref{fig:3ddivfree} as
	\begin{subequations}
		\begin{align}
			%%%% D_x%%%%
			\avg{D}{n+1}{x;i+\thalf, j, k} &= \avg{D}{n}{x;i+\thalf, j, k} - \dt  ~\avg{J}{n+\thalf}{x;i+\thalf, j, k} \notag 
			\\+ \frac{\dt}{\dydz} &\Bigg( \dz~ \avg{H}{n+\thalf}{z;~i+\thalf, j+\thalf, k} - \dz~\avg{H}{n+\thalf}{z;~i+\thalf, j-\thalf, k} 
			+ \dy~ \avg{H}{n+\thalf}{y;i+\thalf, j, k-\thalf} - \dy~\avg{H}{n+\thalf}{y;i+\thalf, j, k+\thalf}\Bigg)  \\
			%%%% D_y%%%%
			\avg{D}{n+1}{y;i, j-\thalf, k} &= \avg{D}{n}{y;i, j-\thalf, k} - \dt  ~\avg{J}{n+\thalf}{y;i,j-\thalf,, k} \notag 
			\\+ \frac{\dt}{\dzdx} &\Bigg( \dx~ \avg{H}{n+\thalf}{x;~i, j-\thalf, k+\thalf,} - \dx~\avg{H}{n+\thalf}{x;i ,j-\thalf, k-\thalf,} 
			+ \dz~ \avg{H}{n+\thalf}{z;i-\thalf, j-\thalf, k} - \dz~\avg{H}{n+\thalf}{z;i+\thalf, j-\thalf, k}\Bigg)  \\
			%%%% D_z%%%%
			\avg{D}{n+1}{z;i, j, k+\thalf} &= \avg{D}{n}{z;i, j, k+\thalf} - \dt  ~\avg{J}{n+\thalf}{z;i,j, k+\thalf} \notag 
			\\+ \frac{\dt}{\dxdy} &\Bigg( \dx~ \avg{H}{n+\thalf}{x;~i, j-\thalf, k+\thalf} - \dx~\avg{H}{n+\thalf}{x;i, j+\thalf, k+\thalf} 
			+ \dy~ \avg{H}{n+\thalf}{y;i+\thalf, j, k+\thalf} - \dy~\avg{H}{n+\thalf}{y;i-\thalf, j, k+\thalf}\Bigg) \\
			%	\end{align}
		% \textnormal{and}
		% 		\begin{align}
			%%%% B_x%%%%
			\avg{B}{n+1}{x;i+\thalf, j, k} &= \avg{B}{n}{x;i+\thalf, j, k} - \dt  ~\avg{M}{n+\thalf}{x;i+\thalf, j, k} \notag 
			\\- \frac{\dt}{\dydz} &\Bigg( \dz~ \avg{E}{n+\thalf}{z;~i+\thalf, j+\thalf, k} - \dz~\avg{E}{n+\thalf}{z;~i+\thalf, j-\thalf, k} 
			+ \dy~ \avg{E}{n+\thalf}{y;i+\thalf, j, k-\thalf} - \dy~\avg{E}{n+\thalf}{y;i+\thalf, j, k+\thalf}\Bigg)  \\
			%%%% B_y%%%%
			\avg{B}{n+1}{y;i, j-\thalf, k} &= \avg{B}{n}{y;i, j-\thalf, k} - \dt  ~\avg{M}{n+\thalf}{y;i,j-\thalf,, k} \notag 
			\\- \frac{\dt}{\dzdx} &\Bigg( \dx~ \avg{E}{n+\thalf}{x;~i, j-\thalf, k+\thalf,} - \dx~\avg{E}{n+\thalf}{x;i ,j-\thalf, k-\thalf,} 
			+ \dz~ \avg{E}{n+\thalf}{z;i-\thalf, j-\thalf, k} - \dz~\avg{E}{n+\thalf}{z;i+\thalf, j-\thalf, k}\Bigg)  \\
			%%%% B_z%%%%
			\avg{B}{n+1}{z;i, j, k+\thalf} &= \avg{B}{n}{z;i, j, k+\thalf} - \dt  ~\avg{M}{n+\thalf}{z;i,j, k+\thalf} \notag 
			\\- \frac{\dt}{\dxdy} &\Bigg( \dx~ \avg{E}{n+\thalf}{x;~i, j-\thalf, k+\thalf} - \dx~\avg{E}{n+\thalf}{x;i, j+\thalf, k+\thalf} 
			+ \dy~ \avg{E}{n+\thalf}{y;i+\thalf, j, k+\thalf} - \dy~\avg{E}{n+\thalf}{y;i-\thalf, j, k+\thalf}\Bigg){.}
		\end{align}
		\label{eqn:dnbupdate}
	\end{subequations}
	The reconstructed values for the electric displacement and magnetic induction, as well as their gradients, form the inputs to the multidimensional GRP. The multidimensional GRP is invoked at each edge. As an output, the GRP gives the time evolution of the resolved state that straddles the edge being considered. From this resolved state, we can evaluate the discrete curl of the electric and magnetic fields along each edge to obtain the globally constraint-preserving update in \Cref{eqn:dnbupdate}. We also have to pay attention, of course, to the source terms for the electric current density and the magnetic current density; these terms are usually stiff and should be handled with a scheme that is unconditionally stable. Furthermore, we want the asymptotic behaviour of the discrete update in \Cref{eqn:dnbupdate} to be such that as $\dt \rightarrow \infty$  the discrete treatment of the source terms gives the same asymptotic result as the differential form of the PDE. Such an unconditional stability is also known as L-stability, and we discuss this in a later section.

	%%%--------------------------------- GRP --------------------------------------------------------------%%%%
	
	\section{Design of a multidimensional GRP solver for Maxwell’s equations and linear hyperbolic partial differential equations in general} \label{sec:multidGRP}
	
	Maxwell’s equations can be written as a system of PDE in the following manner
	\begin{align}
		&\dpd{\con}{t} + \dpd{\fx (\con)}{x} + \dpd{\fy (\con)}{y}  + \dpd{\fz(\con)}{x} = \src (\con),\\
		\intertext{where}
		&\con= \begin{pmatrix*}
			\Dx \\ \Dy \\ \Dz \\ \bx \\ \by \\ \bz
		\end{pmatrix*},
		\qquad	\fx = \begin{pmatrix*}
			0 \\ \imuc{xz}\bx+\imuc{yz}\by + \imuc{zz}\bz \\ -\imuc{xy}\bx-\imuc{yy}\by -\imuc{yz}\bz \\ 0 \\ -\iepsc{xz}\Dx-\iepsc{yz}\Dy -\iepsc{zz}\Dz \\ \iepsc{xy}\Dx+\iepsc{yy}\Dy + \iepsc{yz}\Dz
		\end{pmatrix*},
		\qquad	\fy = \begin{pmatrix*}
			-\imuc{xz}\bx-\imuc{yz}\by - \imuc{zz}\bz \\
			0 \\ \imuc{xx}\bx + \imuc{xy}\by + \imuc{xz}\bz \\ \iepsc{xz}\Dx + \iepsc{yz}\Dy + \iepsc{zz}\Dz \\ 0 \\ -\iepsc{xx}\Dx-\iepsc{xy}\Dy - \iepsc{xz}\Dz
		\end{pmatrix*}, \notag
		\\
		&\fz = \begin{pmatrix*}
			\imuc{xy}\bx+\imuc{yy}\by + \imuc{yz}\bz \\ -\imuc{xx}\bx-\imuc{xy}\by -\imuc{xz}\bz \\ 0 \\ -\iepsc{xy}\Dx-\iepsc{yy}\Dy -\iepsc{yz}\Dz \\ \iepsc{xx}\Dx+\iepsc{xy}\Dy + \iepsc{xz}\Dz \\ 0
		\end{pmatrix*},
		\qquad
		\src = \begin{pmatrix*}
			-\sigma(\iepsc{xx}\Dx + \iepsc{xy}\Dy + \iepsc{xz}\Dz) \\
			-\sigma(\iepsc{xy}\Dx+\iepsc{yy}\Dy +\iepsc{yz}\Dz)\\
			-\sigma(\iepsc{xz}\Dx + \iepsc{yz}\Dy + \iepsc{zz}\Dz)\\
			-\sigma^*(\imuc{xx}\Dx + \imuc{xy}\Dy + \imuc{xz}\Dz) \\
			-\sigma^*(\imuc{xy}\Dx+\imuc{yy}\Dy +\imuc{yz}\Dz)\\
			-\sigma^*(\imuc{xz}\Dx + \imuc{yz}\Dy + \imuc{zz}\Dz)
		\end{pmatrix*}. \notag
	\end{align}
	In the above equations, $\imuc{ij}$ and $\iepsc{ij}$ represents different components of $\imu$ and $\ieps$ tensors where $\imu$ and $\ieps$ are inverses of \num{3x3} symmetric electric permittivity tensor $\bm{\varepsilon}$ and symmetric magnetic permeability tensor $\bm{\mu}$, respectively.
	
	In light of the linearity of the fluxes and source terms in Maxwell’s equations, the above equation can be written in terms of the Jacobians of the fluxes and the Jacobian of the source terms as follows
	\begin{align}
		\dpd{\con}{t} + \ma \dpd{\con}{x} + \mb \dpd{\con}{y}  + \mc \dpd{\con}{z}  = -\sm \con.
		\label{eqn:linearizedhyp}
	\end{align}
	In the above equation, $\ma, ~ \mb ~\textnormal{and} ~ \mc$ are solution-independent  characteristic matrices obtained from the $x, y, z$-fluxes. Likewise, $\sm = -\dpd{\src (\con)}{\con}$ is a constant matrix, where the negative sign has been introduced just to respect the fact that the current terms in Maxwell’s equations are written with a negative sign in front.
	%(Note, while B denotes the magnetic induction vector, B is a characteristic matrix.)
	We would like to design a multidimensional approximate GRP solver for \Cref{eqn:linearizedhyp} specializing it to Maxwell’s equations.
	
	{To describe the development of multidimensional GRP in a step-by-step manner,  we split this section into several parts. In \Cref{subsec:1dgrp}, we briefly describe the $1$D Riemann problem (RP) and generalized Riemann problem (GRP) solvers. This provides us with the lead in to multidimensional Riemann solvers and the multidimensional GRP solver without any source term that we describe in \Cref{subsec:2dgrp}. In \Cref{subsec:GRPwostiff}, we show how this can be used to obtain a GRP without a source term. In \Cref{subsec:GRPwstiff}, we show how the solution of the GRP is obtained in presence of a linear stiff source term.}
	
	\subsection{One dimensional Riemann problem and generalized Riemann problem solvers for linear system}\label{subsec:1dgrp}
	
	A one-dimensional Riemann solver operates at the faces of a mesh because that is where the one-dimensional discontinuities can be diagnosed on a mesh. It takes the two states at a face as input states and provides the resolved state and one-dimensional flux as output.

	\begin{figure}
		\begin{center}
			\includegraphics[width=0.96\textwidth]{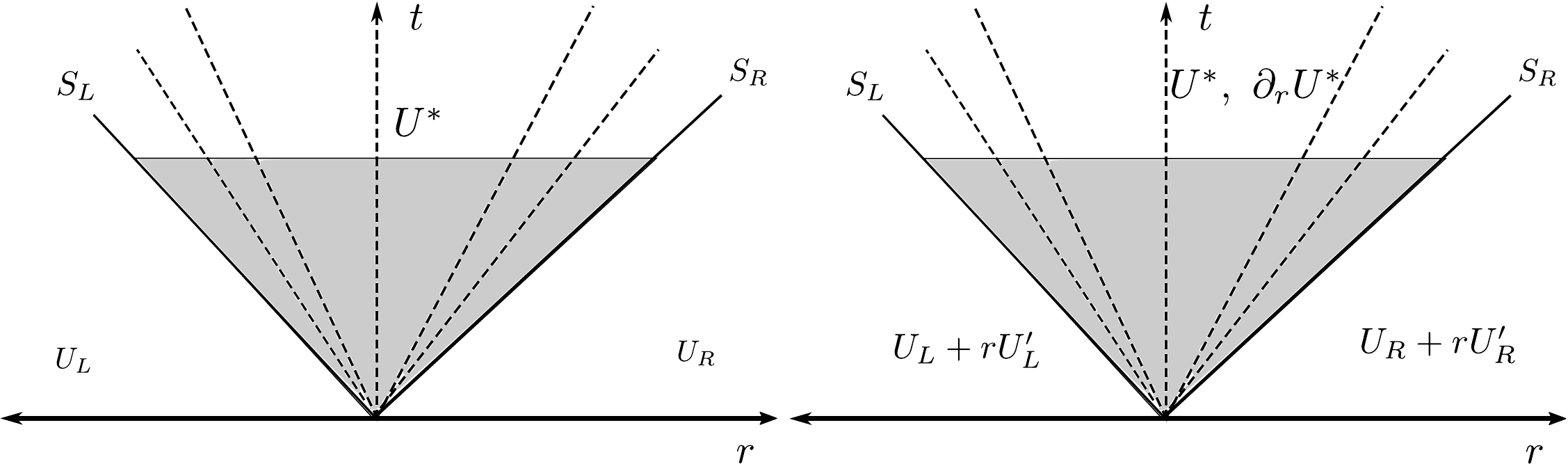}
			\caption{One dimensional solvers that operate at the faces of a mesh. Left panel: One dimensional approximate Riemann solver. It takes two states $U_L$ and $U_R$ at the face as inputs and provides a single resolved state $U^*$ as we use HLL Riemann solver here. Right panel: One dimensional Generalized Riemann solver for linear hyperbolic system. It takes left state $U_L$, right state $U_R$ and their derivatives $\partial_r U_L$ and $\partial_r U_R$ as inputs where $r$ represents any one of the $x,~y,~z$ direction in a Cartesian mesh and provides a resolved state $U^*$ and the derivative $\partial_r U^*$ as outputs. As for linear hyperbolic systems (like Maxwell's equations) maximal wave speeds $S_L,S_R$ are constant, the characteristics curves become straight lines even for GRP.}
			\label{fig:2drpgrpphy}
		\end{center}
	\end{figure}
	
	Analogously, a one-dimensional GRP solver also operates at the faces of a mesh. However, it takes the two states at a face, as well as their spatial gradients, as input states and provides the resolved state and one-dimensional flux and the gradient of the resolved state, as output. The output can then be used to extend the resolved state and its fluxes in time. Please note that we consider the approximate HLL Riemann solver here, which produces only one constant intermediate state between two interacting states. The expression for resolved state is given  by \cite{harten_upstream_1983}:
	{ \begin{align}
			\con^* = -\frac{1}{S_R - S_L}\sbr{(\ma  - S_R \mI)\con_{R} - (\ma -S_L \mI) \con_{L}}, \quad \mI: \textnormal{Identity Matrix},
			\label{eqn:onedhll}
		\end{align}
		where $S_R$ and $S_L$ are extremal speeds in right and left directions.
		
		For a linear hyperbolic system, as $S_R$ and $S_L$ is constant, we can find the derivatives of the resolved state analytically with respect to any arbitrary $r$ direction where $r$ can be any one of $x,~y,$ or $z$ in a Cartesian mesh, and it is given by,
		\begin{align}
			\partial_{r}\con^* = -\frac{1}{S_R - S_L}\sbr{(\ma  - S_R \mI)~\partial_{r}\con_{R} - (\ma -S_L \mI) ~\partial_{r}\con_{L}}.
			\label{eqn:linearhllpartial}
		\end{align}
	}
	We use expressions \eqref{eqn:onedhll} and \eqref{eqn:linearhllpartial} in the next section to obtain resolved states and their derivatives from one dimensional Riemann solvers that are required for a complete description of our multidimensional GRP-based scheme for linear hyperbolic systems.
	
	\subsection{Multidimensional Riemann solver and generalized Riemann solver for linear system} \label{subsec:2dgrp}
	
	A multidimensional Riemann solver operates at the edges of a mesh because that is where the multidimensional discontinuities can be diagnosed on a mesh. We assume a Cartesian mesh to simplify the discussion, but the discussion is indeed generalizable. The multidimensional Riemann solver takes the four states that come together at an edge as input states and provides the resolved state (traditionally called a strongly-interacting state) and multidimensional fluxes as output. Analogously, the multidimensional GRP solver also operates at the edges of the mesh. However, the multidimensional GRP solver takes four states together with their spatial gradients as inputs. As outputs, it produces the strongly-interacting state and multidimensional fluxes, as well as the gradients of the strongly-interacting state. The output can then be used to extend the strongly-interacting state and its fluxes in time.
	
	The edge-based arrangement of electric and magnetic fields for CED in \Cref{fig:3ddivfree} shows that the multidimensional GRP solver provides exactly the desired edge-based data at the very location this data is needed. This highlights the special utility of the multidimensional GRP solver for CED and other involution-constrained applications.

	{The GRP solver is two dimensional, because we would like to invoke it at the edges of the mesh. For illustration, we choose the GRP solver invoked at $z$-edge, and as a result, we focus on the $xy$-plane. However, we will retain derivatives with respect to all three axes in \Cref{eqn:linearizedhyp} as we realize that it might be beneficial to retain the variation in the third direction in fully three-dimensional CED problems.}
	
	For a structured mesh, the specification of the multidimensional Riemann problem at the edges of a Cartesian mesh requires the specification of four input states \cite{balsara_multidimensional_2010,balsara_two-dimensional_2012,balsara_multidimensional_2014-1}. These input states at the initial time are called  $\con_{RU}$ (for right-up), $\con_{LU}$ (for left-up), $\con_{LD}$ (for left-down) and $\con_{RD}$ (for right-down). \Cref{fig:2drp}a shows the input states at the $z$-edge of a mesh, where one is looking down along the $z$-axis. As soon as those input states begin to interact, i.e. at a time that is later than the initial time, four one dimensional Riemann problems get established between the four states.
	Therefore, between the states  $\con_{RU}$ and $\con_{LU}$ an $x$-directional Riemann problem gives rise to the resolved state $\con_U^*$; another $x$-directional Riemann problem between the states $\con_{RD}$  and $\con_{LD}$  gives rise to the resolved state $\con_{D}^*$ ; a $y$-directional Riemann problem between the states $\con_{RU}$  and $\con_{RD}$  gives rise to the resolved state $\con_{R}^*$ ; similarly a $y$-directional Riemann problem between the states  $\con_{LU}$ and $\con_{RD}$ gives rise to the resolved state $\con_{L}^*$. \Cref{fig:2drp}b shows how these resolved states from the one-dimensional Riemann problems are established. When these one-dimensional Riemann problems interact, they form another self-similarly evolving strongly-interacting state $\con^*$ which yields corresponding
	$x$ and $y$ fluxes $\fx^*$ and $\fy^*$. \Cref{fig:2drp}b also shows this strongly-interacting state. While \Cref{fig:2drp}a is in physical space, \Cref{fig:2drp}b is best shown in terms of the wave speeds. Please note that we use the approximate HLL Riemann solver which produces only one constant intermediate state between two interacting states. As a result, the four resolved states $\con _R^*, ~ \con _L^*,~\con _U^*,~\con _D^*$ are constant states without sub-structure. Likewise, the state $\con^*$  has no sub-structure.
	
	Because the characteristic matrices are constant, the extremal speeds in the $x$-direction span $ \xi\in \sbr{S_L, ~S_R } $  and in the $y$-direction span $\psi \in \sbr{S_D, ~S_U}$ are also constant.
	$S_R,~S_L$ are extremal speeds associated with characteristic matrix $\ma$ and $S_U,~S_D$ are extremal speeds associated with characteristic matrix $\mb$. For the case of CED with diagonal permittivity and permeability, we have
	
	\begin{align}\label{eqn:maximalwavespeeds}
		& S_R = \textrm{max}(\sqrt{\imuc{zz}\iepsc{xx}},~ \sqrt{\imuc{yy}\iepsc{zz}}), \quad S_L = -S_R, \notag \\
		& S_U = \textrm{max}(\sqrt{\imuc{xx}\iepsc{zz}},~ \sqrt{\imuc{zz}\iepsc{xx}}), \quad S_D = - S_U.
	\end{align}
	So we see that the extremal wave speeds are very easy to calculate for CED.
	
	\begin{figure}
		\begin{center}
			\includegraphics[width=0.96\textwidth]{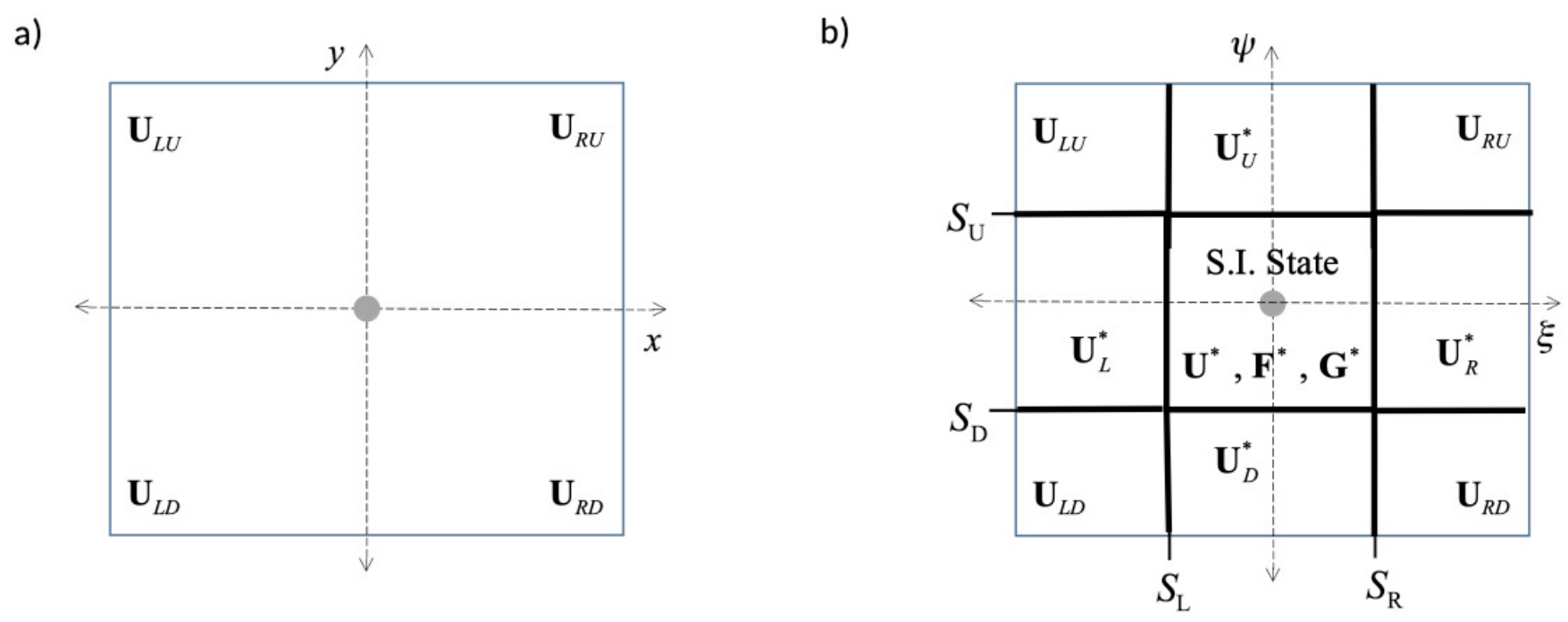}
			\caption{a) Four input states $\con_{LU},~\con_{RU},~\con_{LD},~\con_{RD}$ at the z-edge of a mesh, where one is looking down along the z-axis.
				b) Strongly-interacting state from multidimensional Riemann problem and resolved states from the one-dimensional Riemann problems. While the left panel is in physical space, the right panel is best shown in terms of the wave speeds. Here, {$\xi = x/t$} and {$\psi=y/t$} are the wave speeds in the $x$- and $y$-directions.}
			\label{fig:2drp}
		\end{center}
	\end{figure}

	We now describe how one transitions from a multidimensional Riemann solver to a multidimensional GRP solver. \Cref{fig:2dgrp} shows the input data for the multidimensional GRP – the four input states from \Cref{fig:2drp}a now come in with their gradients in all directions. Therefore, along with the state $\con_{RU}$ we also have its three spatial gradients $\partial_x\con_{RU}, ~\partial_y\con_{RU}, ~\partial_z\con_{RU}$. These gradients can be obtained from any higher order reconstruction {in the neighbouring right-up zone.
		%	{that is used in the right-up zone that abuts the edge being considered}.
		\Cref{fig:2dgrp} shows that similar gradients can be obtained from other neighbouring zones.  We also show the resolved states associated with the one-dimensional Riemann problems and the minimum number of gradients that we should retain in those resolved states. The strongly-interacting state now has all three gradients. Once the gradients have been obtained in the strongly-interacting state, they can be used to obtain the “in-the-small” time-evolution of the strongly-interacting state. This can be done in the sense of a Lax-Wendroff procedure, resulting in a multidimensional GRP solver that is second order in time.
		
		\begin{figure}[htb!]
			\begin{center}
				\includegraphics[width=0.80\textwidth]{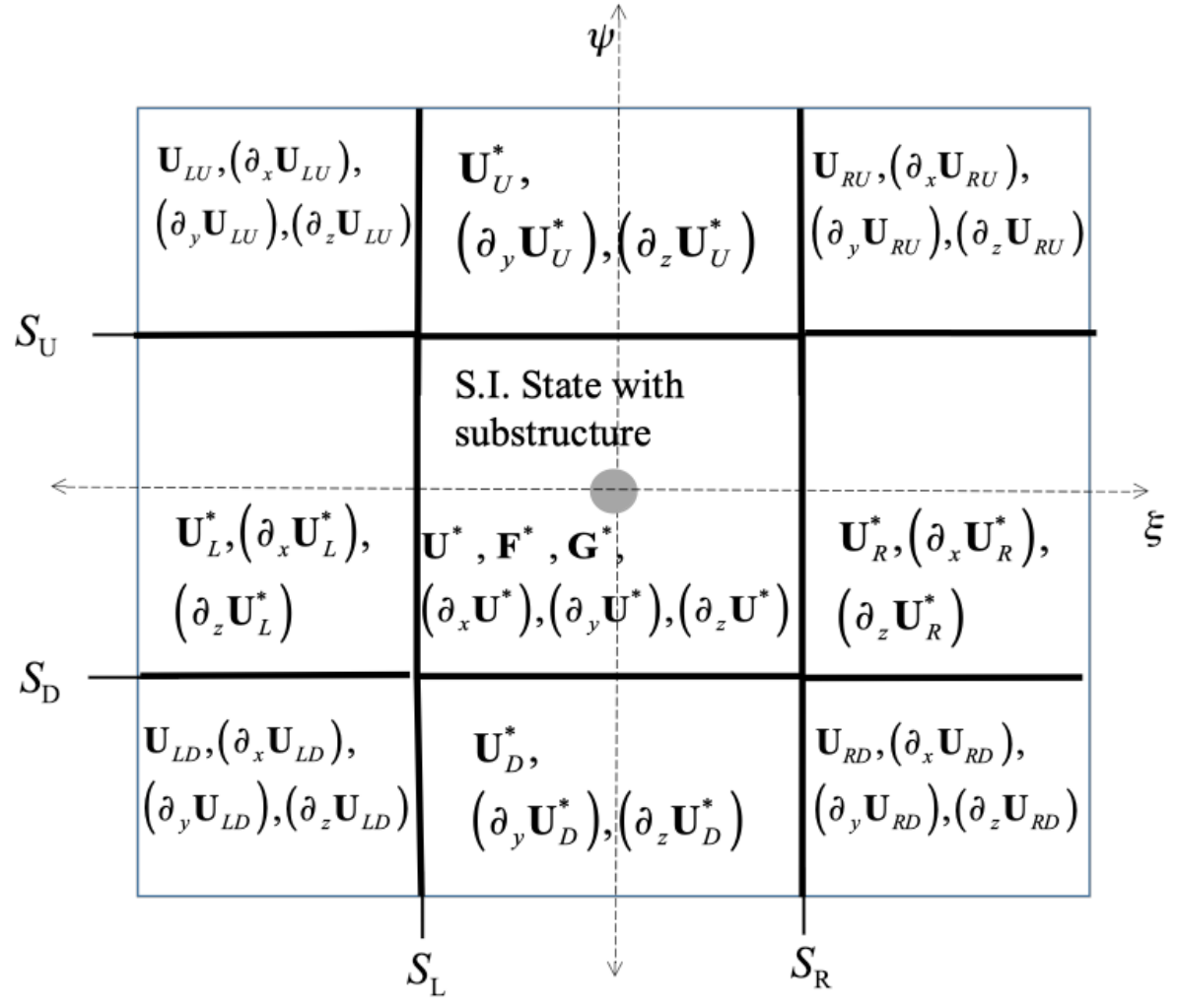}
				\caption{Input data for the multidimensional GRP – the four input states from \Cref{fig:2drp}a now come in with their gradients in all directions. We also show the resolved states of the one-dimensional Riemann problems and the minimum number of gradients that we should retain in those resolved states. The strongly-interacting state now has all three gradients.}
				\label{fig:2dgrp}
			\end{center}
		\end{figure}
		
		First, we focus on the resolved states, represented by $\con_U^*,~\con_D^*,~\con_R^*,~\con_L^*$ in, \Cref{fig:2dgrp} that emerge when we apply the one-dimensional Riemann solvers taking two input states at a time from four input states.  For example, we can obtain the resolved state $\con_U^*$ taking $\con_{RU}$ and $\con_{LU}$ as inputs and using \Cref{eqn:onedhll} and it is given by:
		\begin{align}
			\con_U^* = -\frac{1}{S_R - S_L}\sbr{(\ma  - S_R \mI)\con_{RU} - (\ma -S_L \mI) \con_{LU}}, \quad \mI: \textnormal{Identity Matrix}.
			\label{eqn:resolveduu}
		\end{align}
		Now that we have obtained $\con_U^*$,  we can also find the associated $y$-flux in the upper resolved state as  $\fy_U^* = \mb \con_U^*$.
		
		{{We can see from \Cref{fig:2dgrp}, that we already have $ \partial_y\con_{RU}$ and $\partial_y\con_{LU}$. Also, as we deal with a linear hyperbolic system with constant extremal speeds $S_R$ and $S_L$ here, we can obtain the following expressions for partial derivative of the resolved state in $y$ direction using \Cref{eqn:linearhllpartial}:
				
				\begin{align}
					\partial_y\con_U^* = -\frac{1}{S_R - S_L}\sbr{(\ma  - S_R \mI)~\partial_y\con_{RU} - (\ma -S_L \mI) ~\partial_y\con_{LU}}.
					\label{eqn:difyduu}
				\end{align}
				Similarly, we can obtain the partial derivative of the resolved state with respect to $z$ using \Cref{eqn:linearhllpartial} and it is given by
				\begin{align}
					\partial_z\con_U^* = -\frac{1}{S_R - S_L}\sbr{(\ma  - S_R \mI)~\partial_z\con_{RU} - (\ma -S_L \mI) ~\partial_z\con_{LU}}.
					\label{eqn:difzduu}
				\end{align}
				We can obtain analogous results for $\con_D^*$ by replacing $U$ with $D$ in the subscript of \Crefrange{eqn:resolveduu}{eqn:difzduu}.}}
		
		Now let us focus on the states  $\con_{RU}$ and $\con_{RD}$ and the associated extremal speeds $S_U,~S_D$ in \Cref{fig:2dgrp} and we have a $y$-directional Riemann problem. Using \Cref{eqn:onedhll}, we obtain
		\begin{align}
			\con_R^* = -\frac{1}{S_U - S_D}\sbr{(\mb  - S_U \mI)\con_{RU} - (\mb -S_D \mI) \con_{RD}}.
			\label{eqn:resolvedur}
		\end{align}
		We have obtained $\con_R^*$,  we can also find the associated $x$-flux in the right resolved state as $\fx_R^* = \ma \con_R^*$. The expressions \eqref{eqn:onedhll} and \eqref{eqn:linearhllpartial} are used in the next section to develop a complete multidimensional GRP-based solver.
		
		{{We can see from \Cref{fig:2dgrp} that we already have $ \partial_y\con_{RU}$ and $\partial_y\con_{RD}$ . Therefore, for a linear hyperbolic system and constant speeds $S_U$ and $S_D$, we can obtain $x$ and $z$ derivative of the resolved state $U_R$ using \Cref{eqn:linearhllpartial}.
				\begin{align}
					\partial_x\con_R^* & = -\frac{1}{S_U - S_D}\sbr{(\mb  - S_U \mI)~\partial_x\con_{RU} - (\mb -S_D \mI) ~\partial_x\con_{RD}},
					\label{eqn:difxdur}\\
					\partial_z\con_R^* & = -\frac{1}{S_U - S_D}\sbr{(\mb  - S_U \mI)~\partial_z\con_{RU} - (\mb -S_D \mI) ~\partial_z\con_{RD}}.
					\label{eqn:difzdur}
				\end{align}
				We can obtain analogous results for $\con_L^*$ by replacing $R$ with $L$ in the subscript of \Crefrange{eqn:resolvedur}{eqn:difzdur}.}}
		
		{{At this point, we have all the necessary expressions to obtain the strongly-interacting state $\con^*$ as depicted in \Cref{fig:2dgrp} and it can be obtained using Equation (12) in \cite{balsara_multidimensional_2014-1}}},
		\begin{align}
			\con^* = -\half \Bigg(&\frac{1}{S_R - S_L}\sbr{(\ma  - S_R \mI)\con_{R}^* - (\ma -S_L \mI) \con_{L}^*} \notag \\
			+ &\frac{1}{S_U - S_D}\sbr{(\mb  - S_U \mI)\con_{U}^* - (\mb -S_D \mI) \con_{D}^*}\Bigg).
			\label{eqn:SIstate}
		\end{align}
		{From \Cref{eqn:SIstate}, we can obtain $\fx^* = \ma \con^*$ and $\fy^* = \mb \con^*$ . Also, \Cref{eqn:SIstate} can be formally differentiated in the z-direction to obtain $\partial_z\con^*$ as follows:
			\begin{align}
				\partial_z\con^* = -\half \Bigg(&\frac{1}{S_R - S_L}\sbr{(\ma  - S_R \mI)\partial_z\con_{R}^* - (\ma -S_L \mI) \partial_z\con_{L}^*} \notag \\
				+ &\frac{1}{S_U - S_D}\sbr{(\mb  - S_U \mI)\partial_z\con_{U}^* - (\mb -S_D \mI) \partial_z\con_{D}^*}\Bigg).
				\label{eqn:difzustar}
		\end{align}}
		{For the $x$- and $y$-gradients of the strongly-interacting state, a more sophisticated treatment is described in the ensuing paragraphs. For those who seek the $x$- and $y$-fluxes associated with the state $\con^*$, please see Equations (13) and (14) of \cite{balsara_multidimensional_2014-1}.}
		
		We now focus on introducing $x$- and $y$-gradients in the strongly-interacting state in \Cref{fig:2dgrp}. Consider a general linear hyperbolic system with variation in the $x$-direction. It can be formally written as $$\partial_t \con + \ma \partial_x \con = 0 $$.
		If we differentiate that equation with respect to the $x$-coordinate, it becomes $$\partial_t (\partial_x\con) + \ma \partial_x (\partial_x\con) = 0 $$. We see, therefore, that the $x$-gradient of the solution vector also satisfies a linearized Riemann problem with the same foliation of waves as the original linear hyperbolic system. This insight was first used by Titarev and Toro in \cite{titarev2005ader} to obtain the gradient of the resolved state inside the Riemann fan. It also explains why we will only need $\partial_x \con_R^*$ and $\partial_x \con_L^*$ to obtain $\partial_x \con^*$. Note, however, from an examination of \Cref{eqn:difxdur} that $\partial_x \con_R^*$ and $\partial_x \con_L^*$ do indeed depend on all the $x$-gradients from all the input states. Consequently, we obtain $\partial_x \con^*$  via a genuinely multidimensional contribution from all the input states. Similarly, we will only need $\partial_y \con_R^*$ and $\partial_y \con_L^*$ to obtain $\partial_y \con^*$. In the next two paragraphs, we make this process explicit. We will subsequently provide all possible details using CED as an example.
		
		From the discussion in the previous paragraph, we have understood that the longitudinal (i.e. $x$-directional) gradients of the strongly-interacting state also satisfy the linear system
		
		\begin{align}	\label{eqn:xgradSI}
			&\partial_t (\partial_x\con^*) + \ma \partial_x (\partial_x\con^*) = 0,  \\
			\intertext{with the initial conditions:}
			&\eval[0]{\partial_x\con^*}_{t=0} = \partial_x\con_L^* \quad \textnormal{if}~x<0 \notag \\
			&\eval[0]{\partial_x\con^*}_{t=0} = \partial_x\con_R^* \quad \textnormal{if}~x>0.\notag
		\end{align}
		
		We now use the Titarev-Toro-style linearization. Because the characteristic matrix is constant, the solution of the linear system is easily found. Within the context of the linearization in \Cref{eqn:xgradSI}, we obtain the solution
		\begin{align}\label{eqn:difxustar}
			\partial_x\con^*&=
			\begin{cases}
				\partial_x\con_L^* & \textnormal{when}~ S_L \geq 0\\
				\half\bigg[\partial_x\con_L^*+\partial_x\con_R^*\bigg]+\half \displaystyle\sum_{m=1}^{m_x} \alpha_x^m r_x^m - \half\displaystyle\sum_{m=m_x+1}^{M} \alpha_x^m r_x^m &\textnormal{when}~ S_L <0 <S_R\\
				\partial_x\con_R^* &\textnormal{when}~ S_R \leq 0\\
			\end{cases}
			\\
			\textnormal{with} ~ \alpha_x^m&\equiv l_x^m \cdot \big[\partial_x\con_R^*- \partial_x\con_L^* \big]. \nonumber
		\end{align}
		In the above equation, the eigenvalues $\lambda_x^m,~m=1,2,3, \dots M$ of the left eigenvectors $l_x^m,~m=1,2,3, \dots M$   and the right eigenvectors $r_x^m,~m=1,2,3, \dots M$  are obtained from the characteristic matrix $\ma$ . In \Cref{eqn:difxustar} $m_x$ is defined to be the unique wave for which we have $\lambda_x^{m_x} <0 < \lambda_x^{m_x+1}$. This completes our description of  $(\partial_x\con^*)$. {[We also point out that the omission of a factor of half in front of the eigenvectors in Equation (2.19) of  \cite{balsara_efficient_2018} is indeed an error, and this paper fixes the deficiency in the form of an erratum to that prior paper.]}
		
		Analogous to the discussion in the previous paragraph, the \textit{longitudinal (i.e. $y$-directional) gradients of the strongly-interacting state} also satisfy the linear system
		
		\begin{align}	\label{eqn:ygradSI}
			&\partial_t (\partial_y\con^*) + \mb \partial_y (\partial_y\con^*) = 0,  \\
			\intertext{with the initial conditions:}
			&\eval[0]{\partial_y\con^*}_{t=0} = \partial_y\con_D^* \qquad \textnormal{if}~y<0  \notag \\
			&\eval[0]{\partial_y\con^*}_{t=0} = \partial_y\con_U^* \qquad \textnormal{if}~y>0.\notag
		\end{align}
		As before, we use the Titarev-Toro-style linearization. Because the characteristic matrix is constant, the solution of the linear system is easily found. Within the context of the linearization in \Cref{eqn:ygradSI}, we obtain the solution
		
		\begin{align}\label{eqn:difyustar}
			\partial_y\con^*&=
			\begin{cases}
				\partial_y\con_D^* & \textnormal{when}~ S_D \geq 0\\
				\half\bigg[\partial_y\con_D^*+\partial_y\con_U^*\bigg]+\half \displaystyle\sum_{m=1}^{m_y} \alpha_y^m r_y^m - \half \displaystyle\sum_{m=m_y+1}^{M} \alpha_y^m r_y^m &\textnormal{when}~ S_D <0 <S_U\\
				\partial_y\con_U^* &\textnormal{when}~ S_U \leq 0\\
			\end{cases}
			\\
			\textnormal{with} ~ \alpha_y^m&\equiv l_y^m \cdot \big[\partial_y\con_U^*- \partial_y\con_D^* \big]. \nonumber
		\end{align}
		In the above equation, the eigenvalues $\lambda_y^m,~m=1,2,3, \dots M$ the left eigenvectors $l_y^m,~m=1,2,3, \dots M$   and the right eigenvectors $r_y^m,~m=1,2,3, \dots M$  are obtained from the characteristic matrix $\mb$ . In \Cref{eqn:difyustar} $m_y$ is defined to be the unique wave for which we have $\lambda_y^{m_y} <0 < \lambda_y^{m_y+1}$. This completes our description of  $(\partial_y\con^*)$.
		
		{For CED, the eigenvalues and orthonormal eigenvectors have been documented in \cite{balsara_computational_2017}}. Two of the waves in this system become non-evolutionary because they correspond to the constraints. The case where the permittivity and permeability are diagonal is very important. In that case, we give explicit expressions for $(\partial_x\con^*)$ and $(\partial_y\con^*)$. For $(\partial_x\con^*)$, we have
		
		{\begin{align}\label{eqn:difxustarformula}
				\partial_x\con^*&= \half \Big( \partial_x\con_L^* + \partial_x\con_R^*\Big) + \half  \begin{bmatrix}
					0 \\
					-\sqrt{\dfrac{\tilde{\mu}_{zz}}{\tilde{\varepsilon}_{yy}}}
					\Big((\partial_x\con_R^*)_6  - (\partial_x\con_L^*)_6 \Big) \\
					\sqrt{\dfrac{\tilde{\mu}_{yy}}{\tilde{\varepsilon}_{zz}}}
					\Big((\partial_x\con_R^*)_5  - (\partial_x\con_L^*)_5 \Big) \\
					0\\
					\sqrt{\dfrac{\tilde{\varepsilon}_{zz}}{\tilde{\mu}_{yy}}}
					\Big((\partial_x\con_R^*)_3  - (\partial_x\con_L^*)_3 \Big) \\
					-\sqrt{\dfrac{\tilde{\varepsilon}_{yy}}{\tilde{\mu}_{zz}}}
					\Big((\partial_x\con_R^*)_2  - (\partial_x\con_L^*)_2 \Big) \\
				\end{bmatrix}.
			\end{align}
			In the above equation, $(\partial_x\con^*_{(\cdot)})_i$ denotes the $i$-the component of  the corresponding vector.
			\begin{align}\label{eqn:difyustarformula}
				\partial_y\con^*&= \half \Big( \partial_y\con_D^* + \partial_y\con_U^*\Big) + \half  \begin{bmatrix}
					\sqrt{\dfrac{\tilde{\mu}_{zz}}{\tilde{\varepsilon}_{xx}}}
					\Big((\partial_y\con_U^*)_6  - (\partial_y\con_D^*)_6 \Big) \\
					0\\
					-\sqrt{\dfrac{\tilde{\mu}_{xx}}{\tilde{\varepsilon}_{zz}}}
					\Big((\partial_y\con_U^*)_4  - (\partial_y\con_D^*)_4 \Big) \\
					-\sqrt{\dfrac{\tilde{\varepsilon}_{zz}}{\tilde{\mu}_{xx}}}
					\Big((\partial_y\con_U^*)_3  - (\partial_y\con_D^*)_3 \Big) \\
					0\\
					\sqrt{\dfrac{\tilde{\varepsilon}_{xx}}{\tilde{\mu}_{zz}}}
					\Big((\partial_x\con_U^*)_1  - (\partial_x\con_D^*)_1 \Big) \\
				\end{bmatrix}.
			\end{align}
			The above two equations show us how easy it is to obtain the gradients in the strongly-interacting state.}
		
		\subsection{A GRP solver without any source term}
		\label{subsec:GRPwostiff}
		In its most rudimentary form, a multidimensional GRP solver is used as follows. We hand in the four states and their gradients as inputs to the multidimensional GRP at time $t^n$. The GRP in turn produces the strongly-interacting state $\con^*$ and its gradients $\partial_x\con_U^*, ~\partial_y\con_U^*$ and $\partial_z\con_U^*$ at each $t^n$. However, a GRP solver should enable us to take a temporally second order accurate time update to a time $t^{n+1} = t^n + \dt $ in one call to the GRP solver. As a result, we want the time-centered solution at a time of $t^n + \sfrac{\dt}{2}$. This is obtained by a Lax-Wendroff-like procedure as follows:
		\begin{align}\label{eqn:ustarhalf}
			\con^{*,\shalf} = \con^* - \dfrac{\dt}{2}\Big[ \ma (\partial_x\con^*) +
			\mb (\partial_y\con^*) + \mc (\partial_z\con^*)\Big].
		\end{align}
		With $\con^{*,\shalf}$ in hand, we can easily obtain time-centered electric and magnetic fields at the edges of the mesh. This enables us to find the electric field and magnetic field components at the edges of the mesh with the result that \Cref{eqn:dnbupdate} then be used to make a single-step, second order in time update.
		\subsection{A GRP solver for linear systems with stiff linear source}
		\label{subsec:GRPwstiff}
		The output from the GRP solver will be a state $\con^{*,\shalf}$  updated to a time $t^{n+\shalf}$ at the edges of the mesh, as discussed in the previous sub-section. This state has to be at least first order accurate in time for the overall time update in \Cref{eqn:dnbupdate} to be second order accurate in time. Now let us consider the inclusion of the source term in \Cref{eqn:linearizedhyp}. We can write the update that is analogous to \Cref{eqn:ustarhalf}, but this time we write it formally so that the effect of the source term is included at least up to first order of accuracy. We therefore write
		\begin{align}\label{eqn:ustarhalfws}
			\con^{*,\shalf} = \mathbf{g}(\dt~\sm) \Bigg[\con^* - \dfrac{\dt}{2}\Big[ \ma (\partial_x\con^*) +
			\mb (\partial_y\con^*) + \mc (\partial_z\con^*)\Big]\Bigg].
		\end{align}
		Here $\mathbf{g}(\dt~\sm)$ is a matrix function that depends only on the matrix $\dt ~ \sm$ because the source terms are linear. The matrix function  can consist of any reasonable approximation of $\e^{-(\tfrac{\dt}{2}~\sm)}$. In the next paragraph, we will examine the concept of L-stability in CED. We will then specialize $\mathbf{g}(\dt~\sm)$  to ensure L-stability so that the overall timestep has this very beneficial stability property.
		
		The update in \Cref{eqn:dnbupdate} can be formally written as
		\begin{align}\label{eqn:uupdatenp1}
			\con^{n+1} = \con^n - \dt~\sm ~\con^{*,\shalf} -\dt \mathrm{\mathbf{R}}(\con^{*,\shalf}),
		\end{align}
		where $\con^{*,\shalf}$, $\dt~\sm ~\con^{*,\shalf}$ and  $\mathrm{\mathbf{R}}(\con^{*,\shalf})$ represent the output state from the multidimensional GRP updated to a time $t^n + \sfrac{\dt}{2}$ at the edges of the mesh, the source terms  and a discrete representation of the curl type operator in \Cref{eqn:dnbupdate} respectively. In fairness, the facial currents in \Cref{eqn:dnbupdate} are obtained by averaging the currents provided at the edges of the mesh by the multidimensional GRP. However, that point of detail does not affect the following analysis. The demonstration of L-stability does not rely on the form of the curl-type terms, and so we will ignore the presence of the curl-type terms for the rest of this formal demonstration of L-stability. In other words, to demonstrate L-stability, we will ignore terms with any spatial gradients and focus only on the source terms. When all spatial gradients are set to zero, we have $\con^* = \con^n$. Ignoring all spatial gradients, \Cref{eqn:ustarhalfws,,eqn:uupdatenp1,,} give us
		
		\begin{align}\label{eqn:updateampl}
			\con^{n+1} = \Big[\mI - \dt~\sm ~\mathbf{g}(\dt~\sm)\Big]\con^n = \mathbf{G}(\tilde{\Sigma} \dt) \con^n,
		\end{align}
		where $\mathbf{G(\cdot)}$ is the overall amplification factor of the scheme.
		
		The effect of finite conductivity in Maxwell’s equations is such that, if the spatial gradients do not act, and if the entire system is governed by non-zero conductivity, then the end result after a significantly long time interval should be a zero electric displacement and a zero magnetic induction. In other words, as {$\dt~\sm  \to \infty$, we want $\Big[\mI - \dt~\sm ~\mathbf{g}(\dt~\sm)\Big] \to 0$}. This is a physics-based interpretation of L-stability. The matrix function $\mathbf{g}(\dt~\sm)$ that we choose should reflect that fact. Notice that $\sm$ is either a non-negative diagonal matrix, or it can be diagonalized into such a form via a similarity transformation.
		As a result, we can define $\chi\equiv \dt~\sm$ and simplify our analysis by treating  $\sm$ as a scalar. L-stability is therefore equivalent to demanding that
		
		\begin{align}\label{eq:lstable}
			\lim\limits_{\chi \to \infty} \chi \mathbf{g}(\chi) = 1.
		\end{align}
		In the next three paragraphs we explore different reasonable forms for $\mathbf{g}(\chi)$ so as to finally obtain an L-stable formulation. \Cref{fig:stiffsource} shows the overall amplification factor $\mathbf{G}(\chi) = 1- \chi \mathbf{g}(\chi)$ for different choices.
		
		Let us first examine the exact solution operator of the stiff source term. In that case, and with all gradients set to zero, for \Cref{eqn:ustarhalfws} we obtain the following:
		\begin{align}
			\con^{*,\shalf} = \e^{-(\sfrac{\dt}{2}~\sm)} \con ^n \qquad \textnormal{with} ~\mathbf{g}(\chi) \equiv \e^{-\sfrac{\chi}{2}}.
		\end{align}
		For this choice, we have
		\begin{align}
			\lim\limits_{\chi \to \infty} \chi \mathbf{g}(\chi) = \lim\limits_{\chi \to \infty} \chi \e^{-\sfrac{\chi}{2}}  = 0.
		\end{align}
		So we see that the exact solution operator of the stiff source term does not satisfy the L-stability criteria (\ref{eq:lstable}) and as a result the overall scheme is not L-stable. We might think that the exact evolution operator for the source term should be an ideal choice, but this is not the case when we consider the overall scheme.
		\begin{figure}
			\begin{center}
				\includegraphics[width=0.96\textwidth]{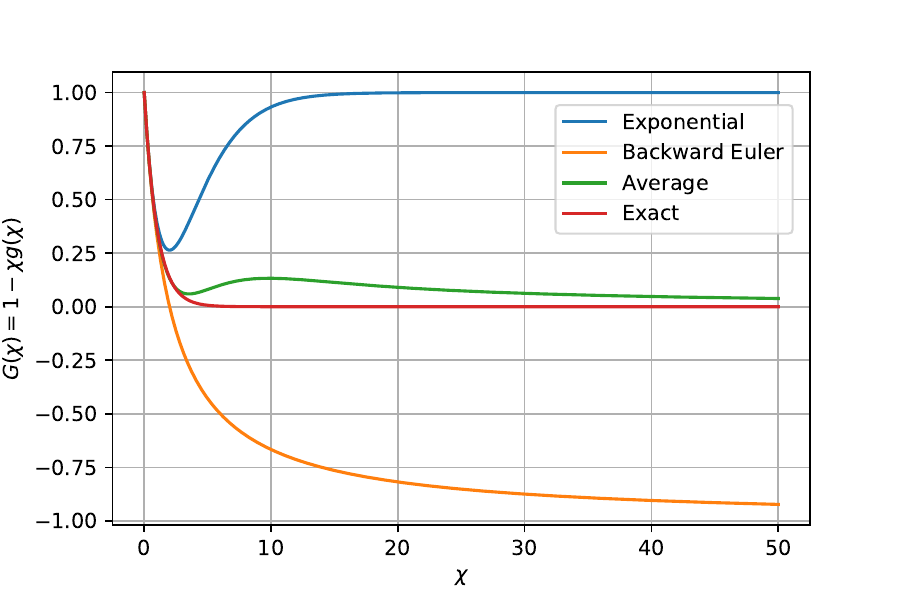}
				\caption{Amplification factor as a function of $\chi = \sm \dt$ . We show that an exact solution as well as the backward Euler solution do not result in an overall time update strategy that has the L-stability property. However, the arithmetic average of the two solutions leads to an overall update strategy that is indeed L-stable. The exact solution is also shown.}
				\label{fig:stiffsource}
			\end{center}
		\end{figure}
		
		Let us next examine the backward Euler solution. In that case, and with all gradients set to zero, we have the following form for \Cref{eqn:ustarhalfws}
		
		\begin{align}
			\con^{*,\shalf} = \Big[\mI + \dfrac{\dt}{2} \sm \Big]^{-1}\con^n \qquad \textnormal{with}~ \mathbf{g}(\chi) \equiv \frac{1}{1+\sfrac{\chi}{2}}.
		\end{align}
		
		For this choice, we have
		\begin{align}
			\lim\limits_{\chi \to \infty} \chi \mathbf{g}(\chi) = \lim\limits_{\chi \to \infty} \frac{\chi}{1+\sfrac{\chi}{2}}  = 2.
		\end{align}
		
		Therefore, we see that the overall amplification factor of the scheme with backward Euler time-stepping $\lim\limits_{\chi \to \infty} \mathbf{G}(\chi) = \lim\limits_{\chi \to \infty} (1- \chi  \mathbf{g}(\chi) ) = -1$.
		
		While the exact and backward Euler solution are not L-stable, the two options begin to hint towards an optimal choice. Since both the exact and the backward Euler solution would give us an overall scheme that is second order accurate  in time, an arithmetic average of the two would also give us a second order accuracy in time.
		If  we take the average of these two options, we obtain the following form for \Cref{eqn:ustarhalfws}
		
		\begin{align}\label{eqn:lstable}
			\con^{*,\shalf} = \half \Big[\e^{-(\dfrac{\dt}{2}~\sm)} +  \Big(\mI + \dfrac{\dt}{2} \sm \Big)^{-1}\Big] \Big[\con^* - \dfrac{\dt}{2} \Big(\ma (\partial_x\con^*) +
			\mb (\partial_y\con^*) + \mc (\partial_z\con^*)\Big)\Big].
		\end{align}
		Now if we set all the gradients to zero, we obtain
		\begin{align}
			\con^{*,\shalf} = \half \Big[\e^{-(\sfrac{\dt}{2}~\sm)} +  \Big(\mI + \dfrac{\dt}{2} \sm \Big)^{-1}\Big] \con^n,
		\end{align}
		and thus
		\begin{align}
			\mathbf{g}(\chi)  = \half \Big(\e^{-(\dfrac{\chi}{2})} + \dfrac{1}{1+\sfrac{\chi}{2}}\Big).
		\end{align}
		For this choice we have
		\begin{align}
			\lim\limits_{\chi \to \infty} \chi\mathbf{g}(\chi)  = \lim\limits_{\chi \to \infty} \frac{\chi}{2} \Big(\e^{-(\dfrac{\chi}{2})} + \dfrac{1}{1+\sfrac{\chi}{2}}\Big) = 1.
		\end{align}
		Therefore, we see that we have found an L-stable scheme and \Cref{eqn:lstable} gives us an expression of the final, successful choice of an overall scheme.

		%%%--------------------------------- Implementation --------------------------------------------------------------%%%%
		
		\section{Pointwise strategy for implementation}\label{sec:Implementation}
		The following steps will result in a one-step, GRP-based, second order accurate in space and time FVTD scheme for CED which preserves the global constraints and is L-stable in the presence of stiff linear source terms. 
		\begin{enumerate}
			\item The primal variables of the scheme are facially-averaged normal components of the electric displacement vector field and the magnetic induction vector field, as shown in \Cref{fig:3ddivfree}.
			These components give us a second order accurate reconstruction of electric displacement vector field and the magnetic induction vector fields following Section III of  \cite{balsara_computational_2017}.
			
			\item Focus on each edge center of a Cartesian mesh. Consider the four zones that abut this edge. Use the reconstructed fields to obtain the four input states to the multidimensional GRP. Because the reconstruction from the previous step also enables us to obtain the gradients in all directions from those four states, we also provide these gradients as inputs to the GRP.
			
			\item Use \Cref{eqn:maximalwavespeeds} to obtain $S_R,~ S_L,~ S_U,~S_D$. This enables us to identify the multidimensional wave model, shown in the left panel of \Cref{fig:2drp}.
			
			\item Use \Cref{eqn:resolveduu,,eqn:difzduu,,eqn:difyduu} to obtain $\con_U^*$  and its gradients in the $y$- and $z$-directions. Do analogously for  $\con_D^*$.
			
			\item Use \Cref{eqn:resolvedur,,eqn:difxdur,,eqn:difzdur} to obtain $\con_R^*$  and its gradients in the $x$- and $z$-directions. Do analogously for  $\con_L^*$.
			
			\item Use \Cref{eqn:SIstate,,eqn:difzustar} to obtain the strongly-interacting state $U^*$ and its gradient in the $z$-direction.
			
			\item Use \Cref{eqn:difxustarformula,,eqn:difyustarformula} to obtain the gradients of the strongly-interacting state  in the $x$  and $y$-directions. Please note that \Cref{eqn:difxustarformula,,eqn:difyustarformula} are just specialized forms of  \Cref{eqn:difxustar,,eqn:difyustar} respectively.
			
			\item If there are no source terms, use \Cref{eqn:ustarhalf}. If source terms are present, use \Cref{eqn:lstable}. This gives us the time-centered states at the edges of the mesh that can be used to construct the curl-type operators in \Cref{eqn:dnbupdate}.
			
			\item If sources are present, obtain the facial current densities by averaging the edge-centered values of the same. This gives us an L-stable treatment of the stiff source terms.
			
			\item Make the update in \Cref{eqn:dnbupdate}. This completes our description of a spatially and temporally second order accurate, globally constraint-preserving, FVTD time update.
		\end{enumerate}
	
	%%%--------------------------------- Accuracy --------------------------------------------------------------%%%%
	
	\section{Accuracy analysis} \label{sec:accuracy}
	\subsection{Propagation of a plane electromagnetic wave in two dimensions}
	\label{subsec:pwave}
	
	In this test problem, we study the propagation of a plane electromagnetic wave through vacuum along the north-east diagonal direction of a two dimensional Cartesian domain spanning $[-0.5, 0.5]\times[-0.5, 0.5]~\si{\square \meter}$. For a detailed description of the problem set up and the electromagnetic field initialization, the readers are referred to \cite{balsara_high-order_2016} for three dimensional version of this test problem and \cite{balsara_computational_2017} for the two dimensional version of this test problem.  Since the analytical solution is known at any space and time, this test problem is very suitable for accuracy analysis. 
	We use a CFL of \num{0.45} and enforce a periodic boundary condition for this problem. The simulation has been run till a time of  \SI{3.5E-9}{\s} second and a uniform mesh has been used in all the runs presented here. \Cref{tab:pwave} shows the accuracy analysis for this test problem. We can see the algorithm meets it designed accuracy for this problem.
	\begin{table}
		\begin{center}
			\begin{tabular}{ccccccccc}
				\hline
				$N_x\times N_y$ & $\|\Dy^h-\Dy\|_{L^1}$ & Ord & $\|\Dy^h-\Dy\|_{L^\infty}$ &Ord \\
				\hline
				$16\times 16$ & 9.8208e-05 & ---&1.5146e-04&---\\
				$32\times 32$ & 2.2130e-05 & 2.15&3.4776e-05&2.12\\
				$64\times 64$ & 5.5153e-06 & 2.00&8.6592e-06&2.01\\
				$128\times 128$ & 1.3850e-06 & 1.99&2.1753e-06&1.99\\
				\hline
				&$\|B_z^h-B_z\|_{L^1}$ & Ord & $\|B_z^h-B_z\|_{L^\infty}$ & Ord\\
				\hline
				$16\times 16$ & 4.9235e-02 & --- &7.8128e-02&---\\
				$32\times 32$ & 1.1492e-02 & 2.10&1.8000e-02&2.12\\
				$64\times 64$ & 2.8693e-03&2.00&4.5064e-03&2.00\\
				$128\times 128$ & 7.2069e-04&1.99&1.1320e-03&1.99\\
				\hline
			\end{tabular}
		\end{center}
		\caption{Accuracy analysis for the second order GRP-WENO scheme for the propagation of an electromagnetic wave in vacuum. A CFL of \num{0.45} was used. The errors and accuracy in the $y$-component of the electric displacement vector and $z$-component of the magnetic induction are shown.}
		\label{tab:pwave}
	\end{table}
	
	\subsection{Compact Gaussian electromagnetic pulse incident on a refractive disk }
	\label{subsec:gauss}
	
	In this two dimensional test problem, we study the propagation of a compact Gaussian electromagnetic pulse that is incident on a refractive disk of refractive index \num{3.0}. The simulation has been performed on a computational domain spanning $[-7.0, 7.0]\times[-7.0, 7.0]~\si{\square \meter}$. The refractive disk of radius \SI{0.75}{\m} is placed at the center of the computational domain. More details about this problem set up and initialization of the compact Gaussian pulse can be found in \cite{balsara_computational_2017}.
	
	This simulation has been run with a CFL of \num{0.45} and continuative boundary condition is enforced for this problem. We stop this simulation at a final time of  \SI{2.33E-08}{\s}. For the simulations presented here, we use a uniform mesh with zones ranging from \num{120 x 120 } to \num{960 x 960}. Since the problem has no analytic solution, we use a \num{1920 x 1920} mesh solution as the reference solution for computing the $L^1$ and $L^\infty$ errors. \Cref{tab:gpulse} shows the result of the accuracy analysis for this problem. The results show that, even for this problem, our algorithm meets its design accuracy.
	\begin{table}
		\begin{center}
			\begin{tabular}{ccccccccc}
				\hline
				$N_x\times N_y$ & $\|\Dy^h-\Dy\|_{L^1}$ & Ord & $\|\Dy^h-\Dy\|_{L^\infty}$ &Ord \\
				\hline
				$120\times 120$ & 4.5549e-05 & --- & 1.7269e-02&---\\
				$240\times 240$ & 2.8758e-05 & 0.66&1.7414e-02&-0.01\\
				$480\times 480$ & 1.1046e-05 & 1.38&6.3582e-03&1.45\\
				$960\times 960$ & 2.3549e-06 & 2.23&1.2758e-03&2.32\\
				\hline
				&$\|B_z^h-B_z\|_{L^1}$ & Ord & $\|B_z^h-B_z\|_{L^\infty}$ & Ord\\
				\hline
				$120\times 120$ & 1.8746e-02&---&2.7283e+00&---\\
				$240\times 240$ & 1.0416e-02&0.85&2.9031e+00&-0.09\\
				$480\times 480$ & 4.1051e-03&1.34&9.6348e-01&1.59\\
				$960\times 960$ & 9.4015e-04&2.13&1.8530e-01&2.38\\
				\hline
			\end{tabular}
		\end{center}
		\caption{ Convergence of error for the second order GRP-WENO scheme for the propagation of a compact Gaussian electromagnetic pulse that is incident on a refractive disk. A CFL of \num{0.45} was used. The errors and accuracy in the $y$-component of the electric displacement vector and $z$-component of the magnetic induction are shown.}
		\label{tab:gpulse}
	\end{table}
	%%%--------------------------------- Numerical Experiments --------------------------------------------------------------%%%%
\section{Test problems} \label{sec:testcases}

\subsection{Refraction of a compact electromagnetic beam by a dielectric slab}\label{subsec:refrac}

In this test problem, we study the refraction of a compact electromagnetic beam impinging on a dielectric slab with a permittivity $2.25 \epsilon_0$ where $\epsilon_0$ is the permittivity of vacuum. Detailed description about this problem set up and the initialization of the electromagnetic beam can be found in  \cite{balsara_computational_2017}.

\sloppy  We perform this simulation on a two dimensional Cartesian domain spanning  ~$[-5.0, 8.0]\times[-2.5, 7.0]~\si{\square \micro \meter}$ using a uniform mesh with $\num{1300x950}$ zones. We use a CFL of \num{0.45 } and stop this simulation at a time of  \SI{4.0E-14}{\s}. The result of the simulation is shown in \Cref{fig:refrac}. The top and bottom rows of \Cref{fig:refrac} shows $B_z,~D_x, ~D_y$ (from left to right) at the initial and final time respectively. The solid vertical black line indicates the interface of the vacuum and the dielectric slab. We have also plotted inclined solid black line to show the angles of incidence, refraction and reflection. These lines are over-plotted with the field components to guide our eye. Since the angle of incidence is \ang{45} for this case, according to Snell’s law, the angle of refraction is \ang{28.12}. We clearly see that our simulation has reproduced the correct value of angle of refraction. 
\begin{figure}
	\begin{center}
		\begin{tabular}{ccc}
			\includegraphics[width=0.33\textwidth]{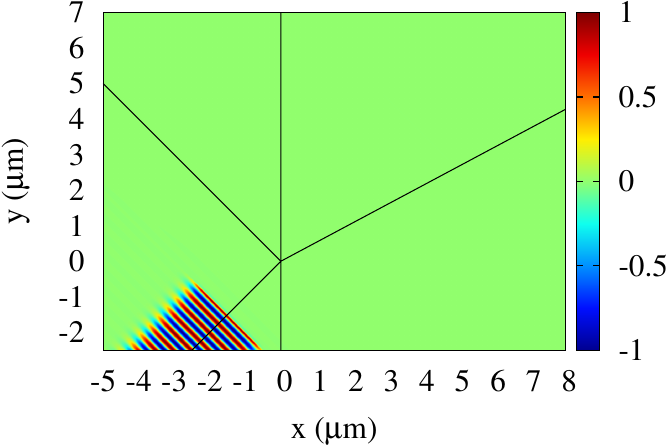} &
			\includegraphics[width=0.33\textwidth]{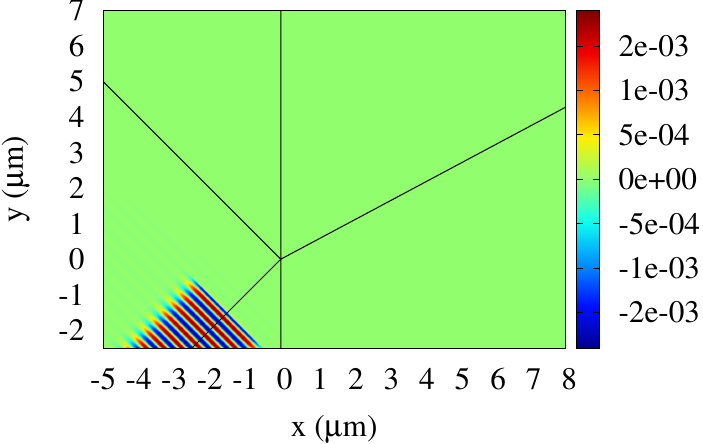} &
			\includegraphics[width=0.33\textwidth]{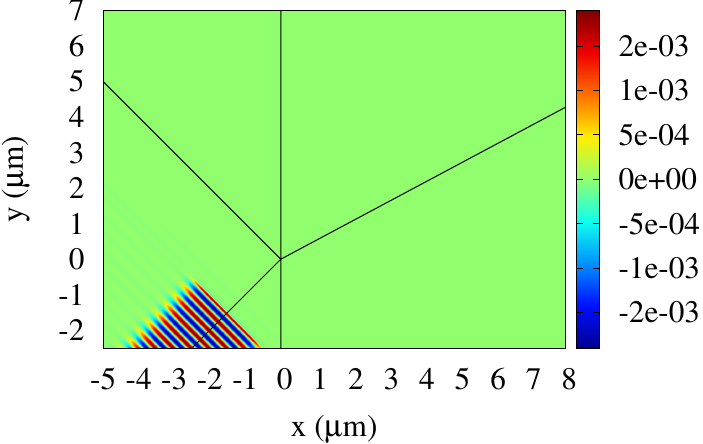} \\
			\includegraphics[width=0.33\textwidth]{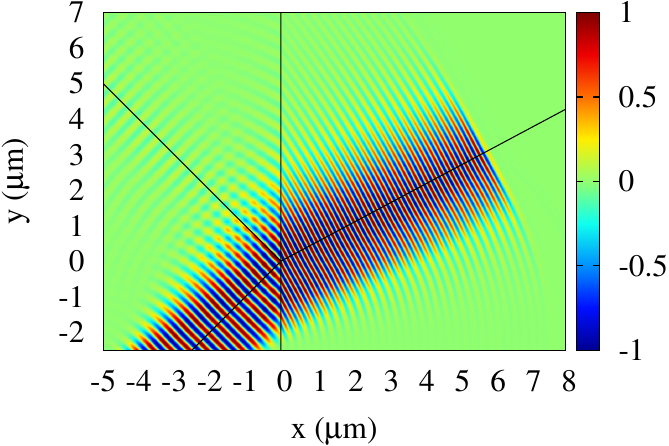} &
			\includegraphics[width=0.33\textwidth]{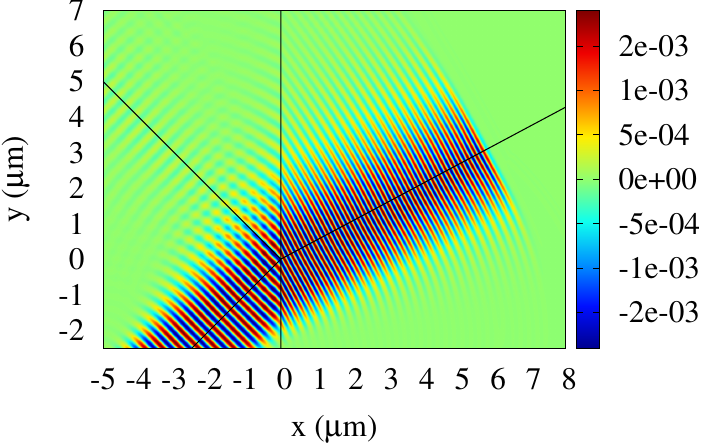} &
			\includegraphics[width=0.33\textwidth]{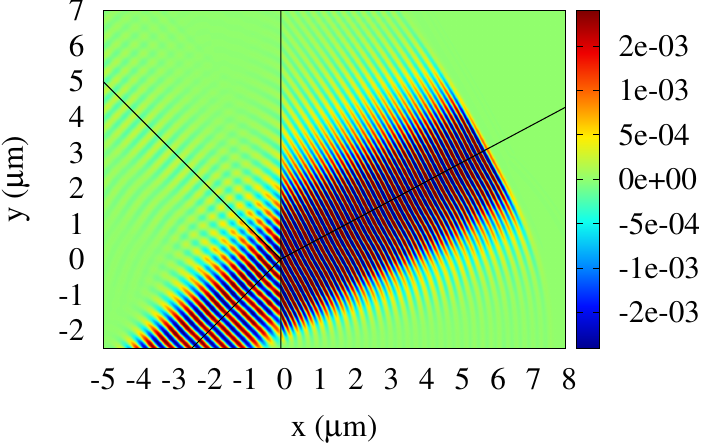} 
		\end{tabular}
		\caption{Refraction of a compact electromagnetic beam by a dielectric slab on a mesh of \num{1300 x 950} cells. The vertical black line indicates the surface of the dielectric slab. The inclined solid black lines demarcate the angle of incidence, the angle of refraction and the angle of reflection. Top row  and bottom row shows $B_z,~ \Dx,~ \textnormal{and}~\Dy$ at the initial time and at final time \SI{4.0E-14}{\s} respectively.}
		\label{fig:refrac}
	\end{center}
\end{figure}

\subsection{Total internal reflection of a compact electromagnetic beam by a dielectric slab}
\label{subsec:reflec}

\begin{figure}
	\begin{center}
		\begin{tabular}{ccc}
			\includegraphics[width=0.33\textwidth]{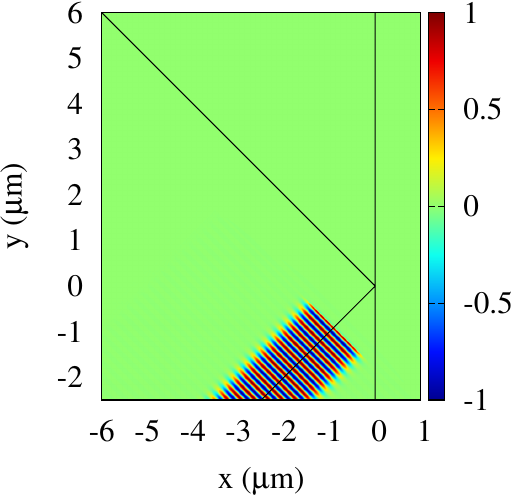} &
			\includegraphics[width=0.33\textwidth]{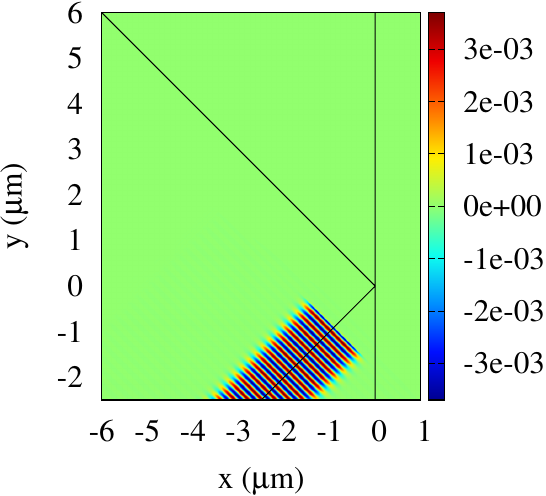} &
			\includegraphics[width=0.33\textwidth]{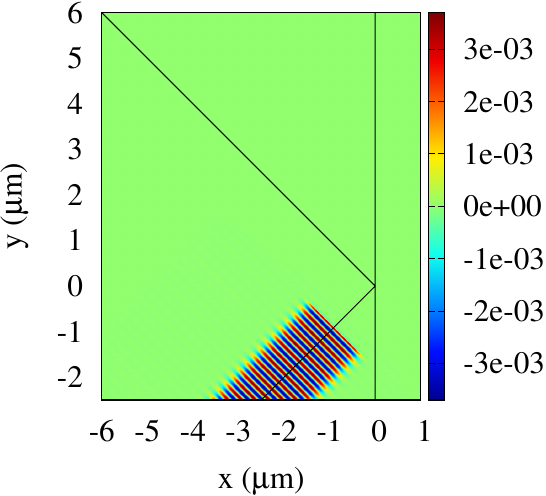} \\
			\includegraphics[width=0.33\textwidth]{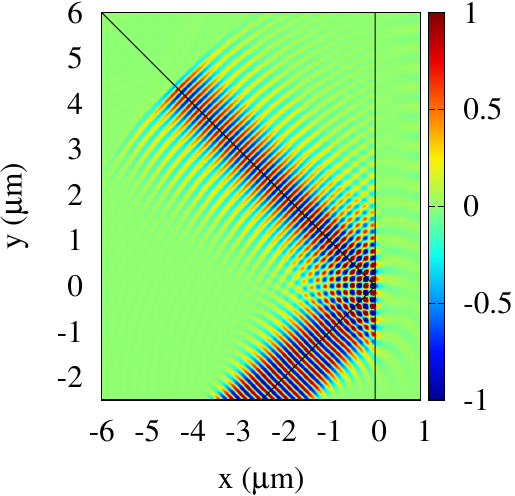} &
			\includegraphics[width=0.33\textwidth]{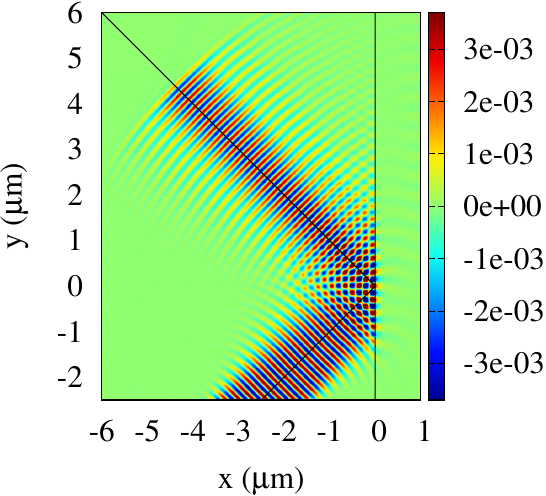} &
			\includegraphics[width=0.33\textwidth]{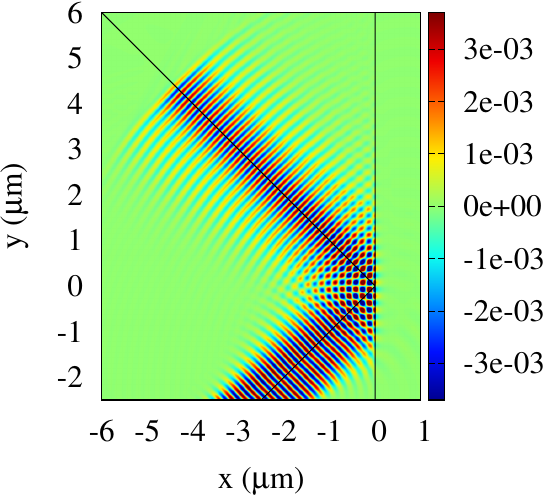} 
		\end{tabular}
		\caption{
			Total internal reflection of a compact electromagnetic beam by a dielectric slab on a mesh of \num{700 x 850} cells: The vertical black line indicates the surface of the dielectric slab. The inclined solid black lines demarcate the angle of incidence and the angle of reflection. Top row  and bottom row shows $B_z,~ \Dx,~ \textnormal{and}~\Dy$ at the initial time and at final time \SI{5.0E-14}{\s} respectively.}
		\label{fig:reflec}
	\end{center}
\end{figure}

In this test problem, we study the total internal reflection of a compact electromagnetic beam when it is incident on the interface separating a dielectric disk of permittivity $4\epsilon_0$ and vacuum at an angle of \ang{45} which is more than the critical angle \ang{30} for such system. For a detailed description of the problem set up, readers are referred to \cite{balsara_computational_2017}.
We perform this simulation on a rectangular xy-domain using a uniform mesh with \num{700 x 850} zones. We use a CFL of \num{0.45} for this simulation run and stop the simulation at a final time of  \SI{5.0E-14}{\s}. The top and bottom rows of \Cref{fig:reflec} show the initial and final configuration of $B_z,~D_x,~D_y$ (from left to right) respectively. The solid vertical black line identifies the interface of the dielectric slab and vacuum. The inclined solid black lines are over-plotted on the field components to guide our eye. The result clearly shows that the incident beam has suffered total internal reflection.

\subsection{Compact electromagnetic beam impinging on a conducting slab}
\label{subsec:elctrobeam}

\begin{figure}
	\begin{center}
		\begin{tabular}{ccc}
			\includegraphics[width=0.31\textwidth]{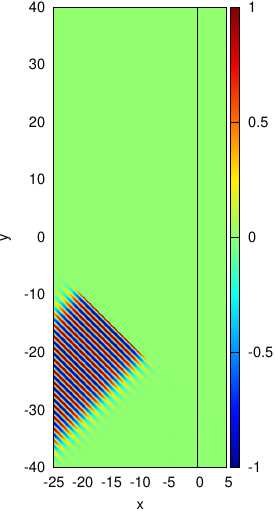} &
			\includegraphics[width=0.33\textwidth]{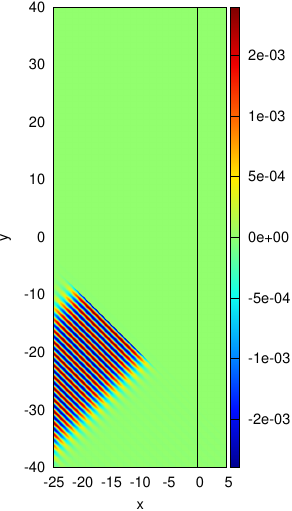} &
			\includegraphics[width=0.33\textwidth]{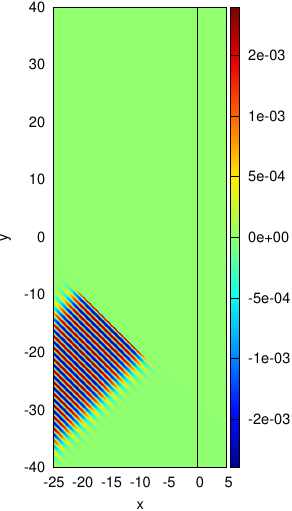} \\
			\includegraphics[width=0.31\textwidth]{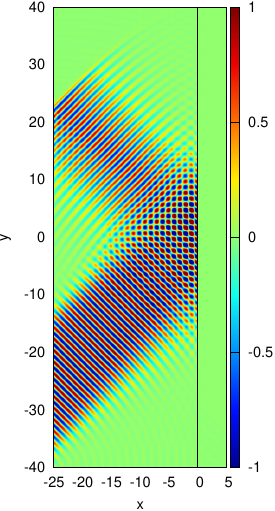} &
			\includegraphics[width=0.33\textwidth]{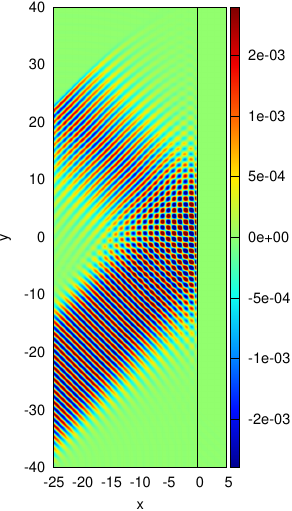} &
			\includegraphics[width=0.33\textwidth]{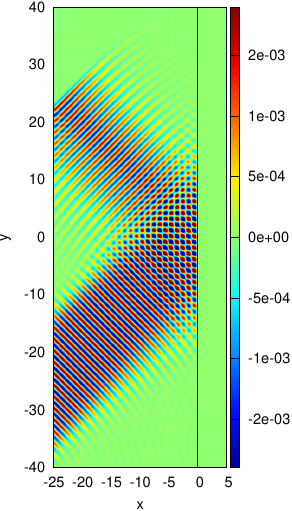} 
		\end{tabular}
		\caption{ Compact electromagnetic beam impinging on a conducting slab. Simulation were performed on \num{1500 x 4000} cells. The conductor is located at $x=0$ in the figure and is shown by the vertical black line. Top row: $B_z,~ \Dx,~ \textnormal{and}~\Dy$ at the initial time. Bottom row: $B_z,~ \Dx,~ \textnormal{and}~\Dy$ at the final time \SI{1.83E-07}{\second} when the beam has reflected off the surface of the conductor. }
		\label{fig:reflecconductor}
	\end{center}
\end{figure}

This test problem is designed to demonstrate the capability of the present algorithm to handle stiff source term. In this test problem, we study the reflection of a compact electromagnetic beam impinging on a slab made up of copper having a conductivity of \SI{5.9E07}{\siemens\per\metre}. The problem set up is described in detail in \cite{balsara_computational_2017}.

We perform this simulation on a rectangular $xy$-domain using a uniform mesh with \num{1500 x 4000} zones. For this simulation, we use a CFL of  \num{0.40} and stop the simulation at a time of  \SI{1.83E-07}{\s}. The top and bottom rows of \Cref{fig:reflecconductor} show the initial and final configuration of $B_z,~D_x,~D_y$ respectively. The solid vertical black line represents the surface of the conducting slab. We can notice the development of interference pattern between the incident wave and the reflected wave close to the surface of the conducting slab.

\subsection{Decay of a sinusoidal wave inside a conductor}
\label{subsec:wavedecay}

Due to the finite skin depth of a conductor, a fraction of the incident wave penetrates it and decays inside it. However, for the meshes used in the previous test problem, we are unable to resolve the skin depth of the copper. Therefore, in this test problem, the simulation set up is designed in a way so that we can resolve the skin depth and study the decay of a sinusoidal wave. The details of the set up and the initialization of the electromagnetic field can be found in \cite{balsara_computational_2017}.

\begin{figure}
	\begin{center}
		\begin{tabular}{cc}
			\includegraphics[width=0.48\textwidth]{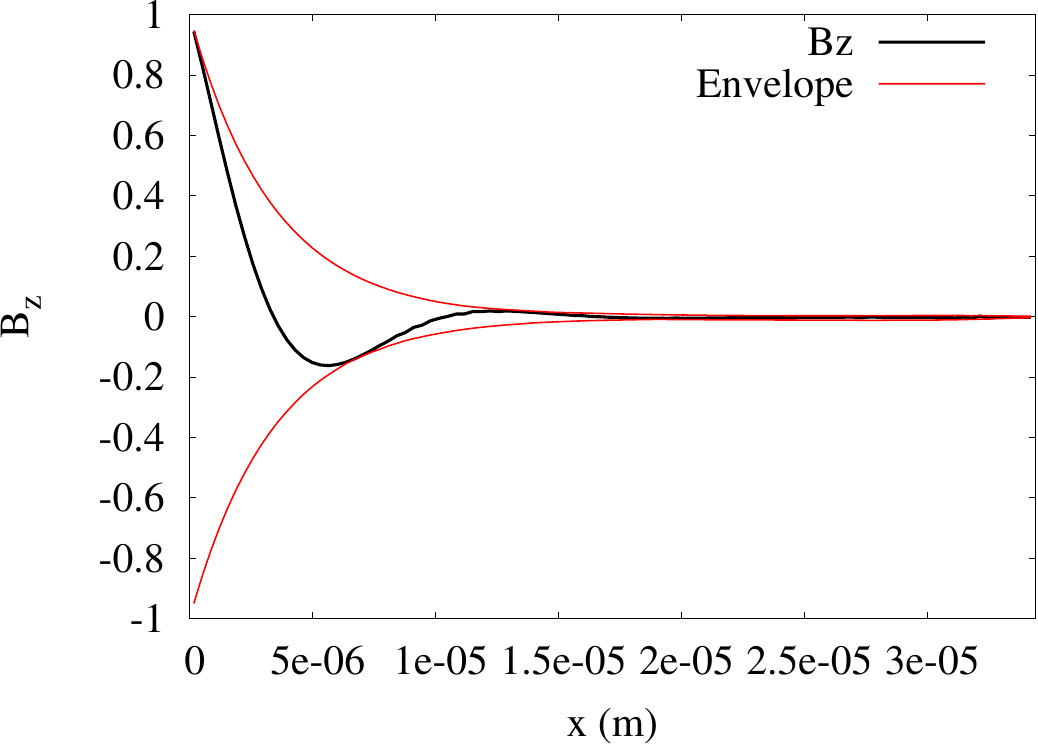} &
			\includegraphics[width=0.48\textwidth]{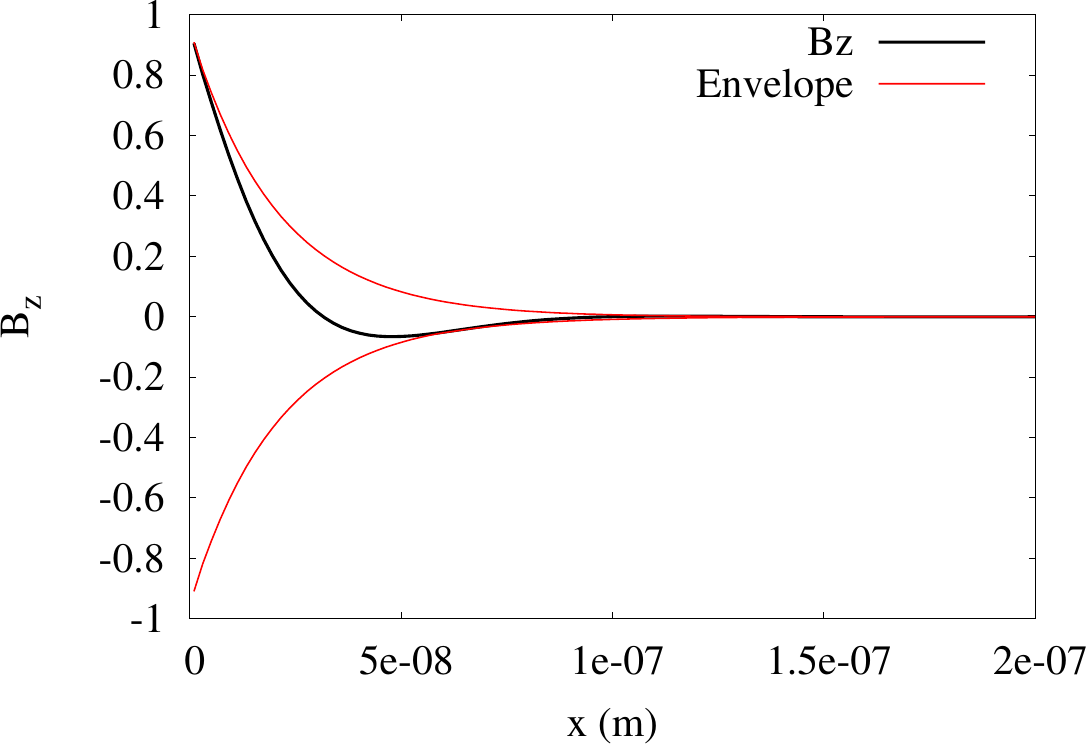} \\
			\includegraphics[width=0.48\textwidth]{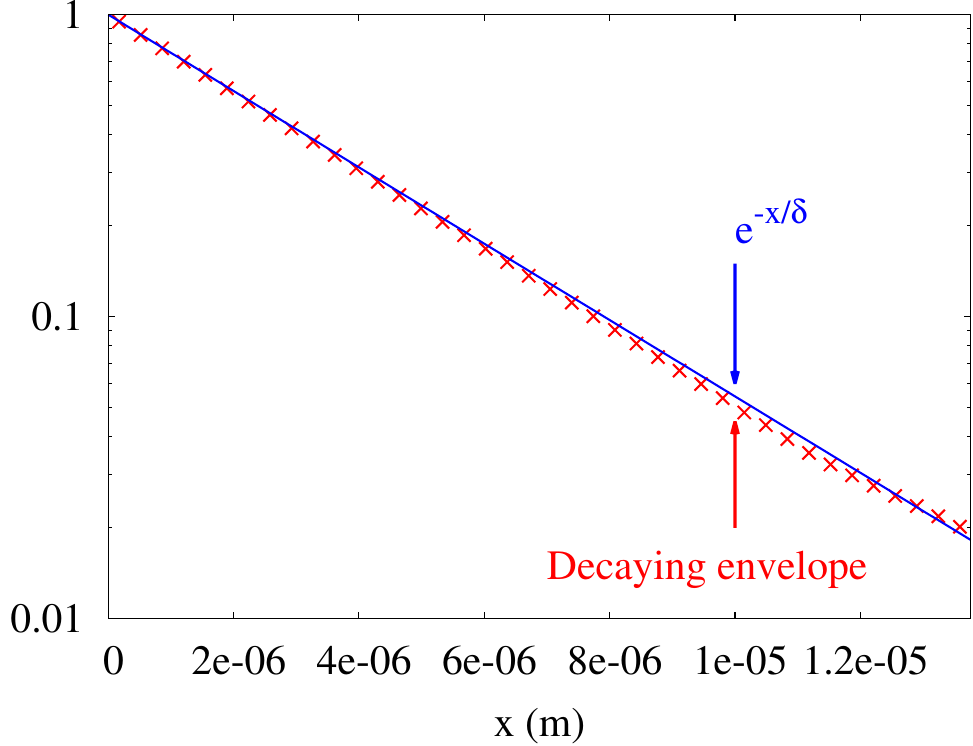} &
			\includegraphics[width=0.48\textwidth]{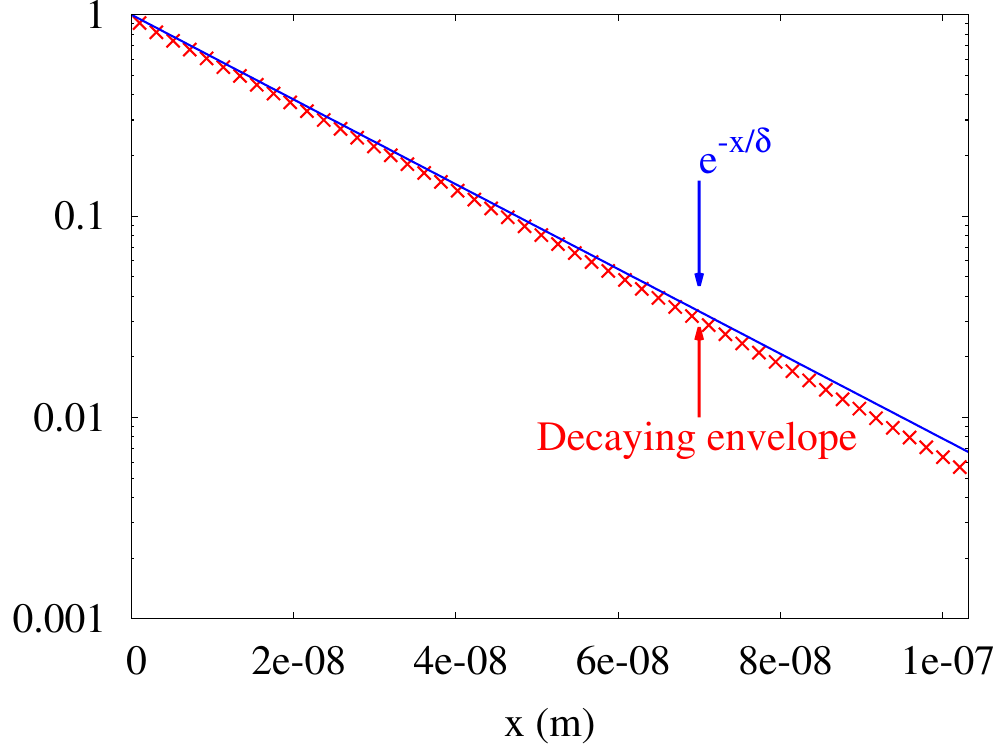}\\
		\end{tabular}
		\caption{ The top left panel and right panel show the radial variations of $B_z$ (black lines) and the decaying
			envelopes (red lines) inside carbon and copper, respectively. The bottom left, and the bottom right panel present the structure of the envelopes (red lines) and the theoretical plots on a semi-log scale for carbon and copper, respectively. }
		\label{fig:decayofwave}
	\end{center}
\end{figure}

Here, we show results of two one dimensional simulations. In one simulation, we study the decay of the sinusoidal wave inside amorphous carbon having a conductivity of  $\sigma = \SI{2.0E03}{\siemens \per \metre} $ and in the other, we study the same inside copper. For both the simulations, we use a one dimensional domain spanning $[0,10\delta]$ where $\delta$ represents the skin depth of the conductor and a uniform mesh with 100 zones. For carbon, we initialize a sinusoidal wave with frequency $\nu = \SI{1.679E13}{\hertz} $, which gives $\delta=\SI{3.44E-06}{\metre}$.  We use a CFL of  \num{0.90} for this run and stop this simulation at a time of \SI{4.76E-13}{\second}. For copper, we initialize a sinusoidal wave with frequency $\nu = \SI{1.0E13}{\hertz} $, which gives $\delta=\SI{2.06E-08}{\metre}$.  We use a CFL of 0.75 for this run and stop this simulation at a time of \SI{4.0E-13}{\second}. Solid black lines in top left and top right panels of \Cref{fig:decayofwave} show the variation of $B_z$ with radial distance inside the carbon and copper, respectively. The solid red lines represent the numerically evaluated decaying envelopes. In bottom left and bottom right panels of \Cref{fig:decayofwave}, we compare the numerically evaluated decaying envelopes (red crosses) with analytically obtained envelopes (blue solid line) on a semi-log scale for carbon and copper, respectively. We can see that our numerical results match very well with the analytical results.  

\subsection{{{Long-distance or long-time propagation of electromagnetic radiation}}}

{{Long-time or long-distance wave propagation is crucial for many problems in electrodynamics. Therefore, it is highly desirable to devise CED schemes with minimal dispersive errors. This test problem
		is designed to demonstrate the numerical dispersion behavior of our numerical scheme. It is also compared with FDTD to prove its superior numerical dispersion behavior over the FDTD scheme. The setup of the problem is analogous to the similar test problem in Section 5.8 of \cite{balsara_computational_2018}. 
		
		To replicate the long-time and long-distance propagation of electromagnetic plane waves, we make electromagnetic plane waves propagate in a small computational domain with periodic boundary conditions in multiple cycles. We choose a computational domain that spans $[-\tfrac{r}{2}, \tfrac{r}{2}] \times [-\tfrac{r}{2}, \tfrac{r}{2}]$ in the $xy$-plane with $r = 6$ divided into a \num{180 x 180} meshes with uniform mesh size and periodic boundary conditions. The exact expression for the electric flux density and the magnetic flux intensity vector fields are given as follows:
		$$\D = c~\epsilon_0\big(-n_y\cos(\phi)\hat{e}_x~+~n_x\cos(\phi)\hat{e}_y\big), \quad \B = \cos(\phi) \hat{e}_z$$ where  $\phi = \frac{2\pi}{n_y}(n_x x + n_y y-ct)$ and $\hat{\bm{n}} = n_x\hat{e}_x + n_y\hat{e}_y$ is the direction of propagation of the plane wave.
		
		As it is well known that wave propagation along the mesh lines or \ang{45} is simpler to replicate, we test the dispersion behaviour of our scheme by choosing the direction of wave propagation along $\hat{\bm n} = (\tfrac{1}{\sqrt{r^2+1}}, \tfrac{r}{\sqrt{r^2+1}})$ which implies that the plane wave are made to propagate at angle $\tan^{-1}(\tfrac{1}{r}) = \tan^{-1}(\tfrac{1}{6}) = \ang{9.462}$ with respect to the $y$-axis.
		
		Simulation was performed until \SI{4.05e-7}{\s} with a CFL number of \num{0.45}. The final time corresponds to \num{20} cycles on the periodic mesh, equivalent to propagating the electromagnetic wave over \num{3600} zones of a uniform mesh.
		
		The left panel of \Cref{fig:dispcomp} shows the variation of $B_z$ normalized by the corresponding amplitude of the sinusoid as a function of $x$ along $y=0$ at the final time. The right panel of \Cref{fig:dispcomp} depicts the same, however, as a function of $y$ along $x=0$. In \Cref{fig:dispcomp}, results obtained using multidimensional GRP, FDTD are compared with analytical solutions, and we can observe that analytical results and multidimensional GRP-computed results are close to each other whereas FDTD-computed results lagged in the left relative to them by \SI{1}{\m} distance. Similar result has also been reported in \cite{balsara_computational_2018}.
		\begin{figure} 
			\begin{center}
				\begin{tabular}{cc}
					\includegraphics[width=0.48\textwidth]{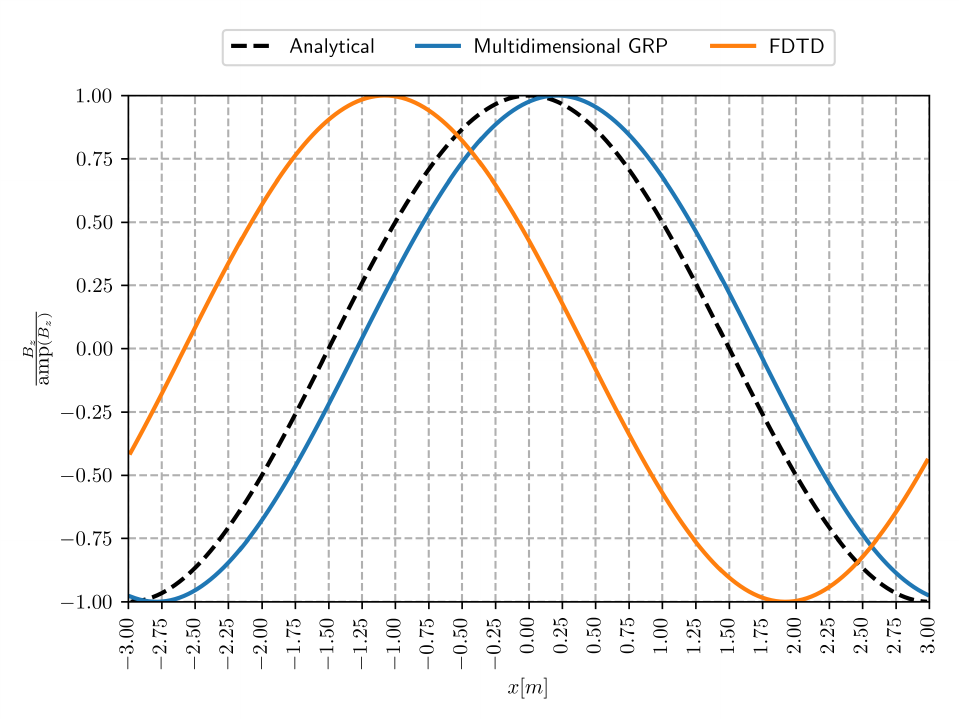} &
					\includegraphics[width=0.48\textwidth]{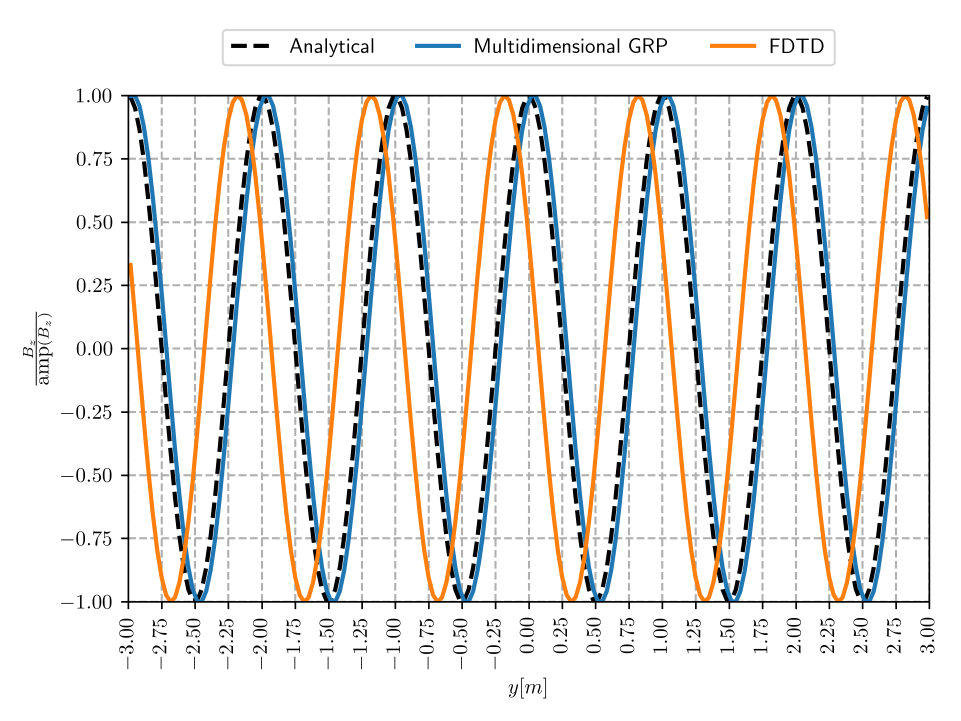}
				\end{tabular}
				\caption{ 
%					Left panel: The variation in $B_z$ normalized by the amplitude as a function of distance along $x$-direction along the line $y = 0$ at the final time for analytical, multidimensional GRP, and FDTD solutions. Right panel: The variation in $B_z$ normalized by the amplitude as a function of distance along $y$-direction along the line $x = 0$ at the final time for analytical, multidimensional GRP, and FDTD solutions.
Long-time propagation of electromagnetic radiation. $B_z$ normalized by the amplitude for multidimensional GRP, FDTD and analytical solution at the final time. Left panel shows the variation along the $x$-axis at $y = 0$, while the right panel shows the variation along the $y$-axis at $x = 0$.}
				\label{fig:dispcomp}
			\end{center}
		\end{figure}
}}

	%%%--------------------------------- Conclusion --------------------------------------------------------------%%%%

\section{Conclusions}\label{sec:conclusion}
In this paper, we have designed an approximate, multidimensional generalized Riemann problem (GRP) solver. The multidimensional Riemann solver takes the four states that come together at an edge as input states and provides the resolved state (traditionally called a strongly-interacting state) and multidimensional fluxes as output. The output can then be used to extend the strongly-interacting state and its fluxes in time. The edge-based arrangement of electric and magnetic fields for CED in \Cref{fig:3ddivfree} shows that the multidimensional GRP solver provides exactly the desired edge-based data at the very location this data is needed. This highlights the special utility of the multidimensional GRP solver for CED and other involution-constrained applications. In this paper, we have designed such an approximate, multidimensional GRP solver for linear hyperbolic systems with stiff, linear source terms. As a result, a one-step update that is temporally second order accurate is achieved.

Our formulation produces an overall constraint-preserving time-stepping strategy based on the GRP that is provably L-stable in the presence of stiff source terms. Our multidimensional GRP formulation, while specialized for CED, is generally applicable to any linear hyperbolic system with stiff, linear source terms.

The multidimensional GRP presented here is intended to be a building block for low dispersion, low dissipation higher order schemes for CED. It could also find utility in aeroacoustics. We also recognize that multidimensional Riemann solvers have found great utility as nodal solvers in Arbitrary Lagrangian Eulerian (ALE) schemes \cite{boscheri_high-order_2014,boscheri_lagrangian_2014}. The availability of multidimensional GRP solvers is expected to greatly simplify the design of higher order ALE schemes because the GRP provides a more accurate trajectory of the strongly-interacting state at each node. Likewise, the Taylor series-based schemes that result from the multidimensional GRP indeed reduce the number of reconstruction stages in ALE schemes. In subsequent papers we will pursue such innovations as well as further develop the field of CED.

%\section*{Declarations}
%
%\subsection*{Funding}
%Dinshaw S. Balsara acknowledges support via NSF grants NSF-ACI-1533850, NSF-DMS-1622457, NSF-ACI-1713765 and NSF-DMS-1821242. Arijit Hazra would like to acknowledge funding support from the Ramanujan Fellowship (RJF/2022/000046) administered by SERB-DST, India.
%
%\subsection*{Competing Interests}
%There is no conflict of interest to declare.
%
%%\subsection*{Author Contributions}
%%Dinshaw S Balsara, Arijit Hazra, Praveen Chandrashekar contributed to the formulation of the problem. Numerical code was developed by Arijit Hazra, Dinshaw S Balsara, Sudip Garain and numerical simulations were performed by Arijit Hazra and Sudip Garain. The primary draft of the manuscript was written by Dinshaw S Balsara and Arijit Hazra. All authors commented on previous versions of the manuscript. All authors read and approved the final manuscript.
%
%\subsection*{Data availability}
%The corresponding author will make numerically simulated data available under reasonable request.
%
%\subsection*{Acknowledgements}
%Arijit Hazra acknowledges support from the Airbus chair on Mathematics of Complex Systems at TIFR-CAM to visit University of Notre Dame. Several simulations were performed on a cluster at UND that is run by the Center for Research Computing.  Computer support on NSF's XSEDE and Blue Waters computing resources is also acknowledged. 

\addcontentsline{toc}{section}{References}
%\addcontentsline{toc}{section}{References}
\bibliographystyle{unsrtnat}
\bibliography{mybibfile}

\end{document}